%% file: ClasTroo.tex
\documentclass[11pt,reqno]{amsart}

\usepackage{amssymb,verbatim,amscd}

\usepackage[cp866]{inputenc}\usepackage[english,russian]{babel}
\usepackage[mathcal]{euscript} 
\let\sM\smallmatrix  \let\endsM\endsmallmatrix
\def\ig{\vphantom{\Bigr)}}  

\usepackage{TRII}
\usepackage{that47}

\title{}

\address{M.M.Graev, Scientific Research Institute for System Studies of Russian Academy of Sciences, Nakhimovsky prosp., 36, korpus 1, Moscow, 117218,  Russia}
\email{mmgraev@gmail.com}

\begin{document}

\def\contentsname{Contents}
\def\refname{References}
\def\bibname{References}
\def\figurename{Fig.}
\def\tablename{Table}
\def\listtablename{List of Tables}

\begin{center}{}\Large\bf
$T$-root systems of exceptional flag manifolds.\\
Classification and invariants
\end{center}


\begin{center}{}\large
M. M. Graev
\end{center}

\footnotetext{Supported by RFBR, grant 13-01-00190. 
}

\maketitle

\vskip-1cm
\vskip-1cm




\begin{quote}
{\small 
\noindent
\textsc{Abstract}.
\let\bar\overline
\input abstract
}
\end{quote}

\begin{quote}
\tableofcontents
\end{quote}
\vskip-1cm
\vskip-1cm

\listoftables
\clearpage

\section*{Introduction} 

\input intro-dv.tex

\null\clearpage

\input T-2

\input T-1
\clearpage

\begin{table} 
\caption{Classification and invariants of T-root systems of ranks $\ge2$}
\small
\textbf{Notations. }
\label{table:31}
{1) For flag manifolds:}
\begin{itemize}
\item $f :$ the number of the standard complex forms of a flag manifold;
\item $[1,1],\dots,[0,0] :$ for flag manifolds of $E_\ell$, $\ell=7,8$;\enskip
we write $[*,1]$ (respectively $[1,*]$) if the semisimple part of the isotropy subgroup is conjugate to a subgroup of the standard $A_{\ell-1}$ (respectively $E_{\ell-1}$).
\end{itemize}
{2) For any $T$-root system $\Omega :$} 
\begin{itemize}
\item $d=\tfrac12\,|\Omega |=|\Omega ^+| $
\item $c=\prod_{\omega  \in \Omega ^+} g_{\omega }$,\quad where
\quad  $g_{\omega }=\max\{k\in\ZZ: \tfrac1k\,\omega \in \Omega \}$
\item $v=|\Omega ^+_V| $,\quad where
\quad $\Omega ^+_V:= \{\omega \,\in\, \Omega ^+ :
2\omega - \gamma  \notin \Omega^+, \,\,\forall\,\gamma \in \Omega \cup (0)
,\gamma \ne \omega\}$. 
\item $w=\prod_{\omega \, \in\, \Omega ^+_V} g_{\omega }$
\end{itemize} 
\end{table}
\input Tables/Tables31

\input Tables/Nmbering
\clearpage


\input Tables/Tables32
\clearpage


\end{document}

%% file: abstract.tex
We classify the systems of $T$-roots of the flag manifolds $M$ of the exceptional compact simple Lie groups with the second Betti number $b_2(M)\ge2$.
 

%% file: intro-dv.tex
\input hat
\NNULL{ }

A flag manifold of a compact semisimple Lie group $G$ is an adjoint orbit  $M = Ad(G).x = G/H$  of $G$. It carries  many invariant   structures, for example, invariant  complex  structures and     associated invariant K\"ahler-Einstein structures, a big family of invariant K\"ahler structures, many non-K\"ahler  invariant Einstein metrics etc. 
An important invariant of a flag manifold is  the  system of $T$-roots.
In the case of the full flag  manifold $M = G/T^{\ell}$ (i.e., a regular orbit of $G$), the system of $T$-roots reduces to the system of roots (w.r.t. a maximal torus $T^{\ell} \subset G$).
In the general case the $T$-root system can be defined as follows.

\medskip

Let $G$ be a compact semisimple Lie group, and $M=G/H$ a flag manifold. The $T$-roots \cite{AP} of $G/H$ are the characters of the torus $T^k=Z(H)$ corresponding to the one-dimensional $T^k$-invariant subspaces of $\CC \otimes (\mathfrak g/\mathfrak h)$. 
The system $\Omega$ of the $T$-roots is a subset of the vector space $\RR^k=(Lie\,(Z(H)))^*$. It is easy that $\Omega$ is not contained in a hyperplane.  The number $k$ is called rank of $\Omega$. We consider $\Omega$ up to linear isomorphism. 

\medskip

In \cite{2006}, see also \S\,\ref{sect:1}, we enumerate the $T$-root systems (up to linear isomorphisms) of the flag manifolds $M$ of the classical simple Lie groups $G$. 
The aim of this paper is to classify the systems of $T$-roots of the flag manifolds $M$ of the exceptional compact simple Lie groups with the second Betti number $b_2(M)\ge2$, that is, with $k\ge2$.

\medskip

Note that every finite root system is a system of $T$-roots. 
E.g., the root system of type $G_2$ is the $T$-root system of flag manifolds $G_2/T^2$, $E_6/T^2\cdot (A_2)^2$, $E_7/T^2\cdot A_5''$, $E_8/T^2\cdot E_6$, and $F_4/T^2\cdot \widetilde A_2$. This is used in the recent paper \cite{SACos} to enumerate positive definite invariant Einstein metrics up to isometry and scale on these manifolds ($3, 7, 7, 7, 7$ metrics respectively, cf. also \cite{11-11}). 

\medskip

The article is organized as follows: 

In \S\,\ref{sect:1} and \ref{sect:2} we state basic definitions and collect some preliminary results.

In \S\,\ref{sect:trootwo} we enumerate irreducible $T$-root systems of rank two (Table~\ref{table:trootwo-Omega}). and the corresponding flag manifolds of Lie groups $E_6$, $E_7$, $E_8$, and $F_4$  (Tables~\ref{table:trootwo-flag} and~\ref{table:trootwo-flag-complex}). We prove consequences for the case of rank $k\ge 3$.

In \S\,\ref{sect:4} we state main classification results and prove some concequences. In particular, we enumerate flag manifolds of $E_6$, $E_7$, $E_8$, and $F_4$ with isomorphic $T$-root systems of ranks $k\ge3$ (Theorem~\ref{THM:FOm}, cf. Table~\ref{table:Numberings}), and describe the $T$-root system with $d\le 10$ positive $T$-roots, corresponding to flag manifolds of $E_6$, $E_7$, $E_8$, $F_4$.

In \S\,\ref{sect:tables} we introduce some simple invariants of $T$-root system and get classification tables~\ref{table:31} and~\ref{table:Numberings}.

\ENDINPUT

%% file: hat.tex
\let\ENDINPUT 
\def\NNULL#1{}

%% file: T-2.tex
\input hat
\NNULL{
\theoremstyle{plain}
\newtheorem*{CLAIM*}{Claim}
}

\section{Definitions. Classification of T-root systems of classical types} \label{sect:1}

Let $G$ be a compact connected semisimple Lie group of rank $\ell$.
We fix a maximal torus $T^\ell \subset G$.

A flag manifold $M=G/H$ of $G$ is the factor space by the centralizer $H$ of some torus in $G$.
Any flag manifold is simply connected.
We will suppose later that $T^{\ell}\subset H$.  In other words, the basepoint $eH\in G/H$ belongs to the finite set of the fixed points of $T^{\ell}$.


More precisely (although to fix the idea we identify $M$ with $G/H$), we will consider $M$ as a manifold with transitive and locally effective action of the group $G$ by ignoring the basepoint. In other words, coset spaces $G/H$ and $G/H_1$ represent the same flag manifold $M$ if and only if subgroups $H$ and $H_1$ of $G$ are conjugate, $H_1=xHx^{-1}$ for some $x\in G$.

\medskip

We associate now with everyone of considered coset spaces $G/H$  a finite vector set, called {\bf T-root system} (\cite{AP}). With a flag manifold $M$ one can associate an equivalence class of such systems (up to linear isomorphisms). Later in this section we enumerate $T$-root systems (up to isomorphisms) of the flag manifolds of the classical groups (by \cite{2006}).

Let $Q_G$ be the root lattice of $\mathfrak g $ (with respect to $T^{\ell}$), and $Q_H\subset Q_G$ the sublattice generated by the roots of $[\mathfrak h,\mathfrak h]$. By definition of the flag manifold $G/H$, the quotient group $Q_G/Q_H$ is torsion-free. More precisely, $Q_G/Q_H$ can be considered as a vector lattice in $\RR^k$, the dual space of the Lie algebra of the torus $T^k=Z(H) :$
$$
\RR^k = \RR\otimes(Q_G/Q_H)=(Lie\,T^k)^*.
$$
\begin{DEF*} 

(\cite{AP}) 
Let $\alpha$ be a root of $\mathfrak g$, but not a root of $[\mathfrak h,\mathfrak h]$, that is, $\alpha\notin Q_H$. Then the natural projection of $\alpha$ into $\RR^k\setminus (0)$ (that is, $\alpha+Q_H$) is called a {\bf T-root}.

\end{DEF*}

E.g., $T$-roots of $G/T^{\ell}$ are roots of the Lie algebra $\mathfrak g$ with respect to $T^{\ell}$.

\smallskip

We call a {\bf T-root system} of $G/H$, or simply a T-root system, the set $\Omega=\Omega_{G/H}$ of all $T$-roots.
The {\bf rank} of $\Omega$ is the number 
$$
k=\dim(Z(H))=b_2(M),
\quad \mbox{the second Betti number of $M$}.
$$

There is a natural correspondence between $T$-roots $\omega \in \Omega $ and irreducible submodules $\mathfrak m_{\omega}$ of the complex $H$-module $\CC\otimes(\mathfrak{g/h}):$ 
$$
\mathfrak m_{\omega} = \bigoplus_{\alpha\,\in\,\omega}\mathfrak g^{\alpha
}.
$$
Indeed, every $\mathfrak m_{\omega}$ is irreducible by the profound Lemma~3.9 in \cite[Chap. 3]{ViGOn}. Moreover, submodules $\mathfrak m_{\omega}$, $\omega \in \Omega$, are pairwise non-equivalent even as $Z(H)$-modules. 
Since $\mathfrak m_{\omega}$ and  $\mathfrak m_{-\omega}$ are complex conjugate, the \textit{real} isotropy $H$-module $\mathfrak{g/h}$ decomposes as a direct sum of
$$
d=\frac12|\Omega|
$$
mutually non-equivalent irreducible submodules.

\begin{DEF*}

A {\bf T-root system in $\RR^k$} is a subset $A\Omega_{G/H}$ for any $A\in\mathrm{GL}(k,\RR)$.
Let $\opn{Aut}(\Omega):=\{A\in\mathrm{GL}(k,\RR): A\Omega=\Omega \}$.  

\end{DEF*}

\begin{DEF*} 

Two $T$-root systems $\Omega$ and $\Omega'$ in $\RR^n$ are \textbf{isomorphic}, if $\Omega=A\Omega'$ for some $A\in\mathrm{GL}(n,\RR)$.

\end{DEF*}

Here is an example of such isomorphic systems.
If two subgroups $H_1=H$ and $H_2$ of $G$ are conjugate, i.e., $H_2=xHx^{-1}$ for some $x\in G$, then there is a commutative diagram
$$
\begin{CD}
Z(H_1) @>{\sim}>> Z(H_2)\\
@VVV  @VVV\\
\mathrm{GL}(\mathfrak g/\mathfrak h_1)@>{\sim}>>\mathrm{GL}(\mathfrak g/\mathfrak h_2)
\end{CD}
$$
Since $\Omega_{G/H}$ can be defined independently of $T^\ell$, starting from the natural action of $Z(H)$ on $\mathfrak g/\mathfrak h$,
this means that $T$-root systems $\Omega_{G/H_i}$, $i=1,2$, are isomorphic. 
 
\smallskip

This observation leads to the following definition.

\begin{DEF*}

We call a {\bf T-root system of the flag manifold} $M$ the isomorphism class $\Omega_M$ of the $T$-root system of $G/H$. We will consider later $\Omega_{G/H}$ up to isomorphism and write $\Omega_M=\Omega_{G/H}$. 

\end{DEF*}

We enumerate in  \cite{2006} the $T$-root systems of the flag manifolds of the classical groups $G$. There are the roots systems
of the types $A_n,B_n,C_n,D_n,BC_n$, and the systems of vectors obtained from $C_n$ and $BC_n$ by removing a part of long roots, as the next theorem states. 

\begin{THM}\label{THM:classtypes} 

The isomorphism classes of the $T$-root systems of the flag manifolds $M$ of the classical compact simple Lie groups are the following:
\begin{multline*}
A_1, \quad BC_1,\quad A_2, \quad B_2,\quad BC_{2,1}, \quad  BC_2,\quad
\\ A_3,\quad C_{3,1},\quad C_{3,2},\quad C_{3},\quad B_{3},\quad BC_{3,1},\quad BC_{3,2},\quad BC_{3},\qquad
\\ A_4,\ldots ,A_n, \quad
D_n, \quad C_{n,1},\ldots ,C_{n,n-1}, \quad  C_n, \quad  \\
B_n, \quad  BC_{n,1},\ldots ,BC_{n,n-1},\quad  BC_n, \quad  A_{n+1}, \ldots,
\end{multline*}
where $C_{n,n-k}$ and $BC_{n,n-k}$ respectively are obtained from $C_n$ and $BC_n$ by removing a $k$ pairs of opposite roots of the maximum length.
Each of these $T$-root systems except 
$D_n$, $n>3$, corresponds to infinite set of flag manifolds $M$, indicated in \cite[Introduction]{2006}.

\end{THM}

\section{Preliminary}\label{sect:2}

\subsection{Standard complex forms of flag manifold}\label{sect:2.1}

Assume that $G^\CC$ is the complex semisimple Lie group with the maximal compact subgroup $G$.
(Here without loss of generality we may assume that $G=Ad(G)$, $G^\CC=Ad(G^\CC)$).
Fix the standard positive Borel subgroup $B_+ \subset G^{\CC}$, so that $B_+\cap G= T^{\ell}$. 

Let $P$ be a parabolic subgroup of $G^{\CC}$ such that
$$
B_+ \subset P \subset G^{\CC}.
$$ 
Remark that if subgroups $P$ and $P'$ are conjugate, and $B_+\subset P'$, then $P'=P$ (this is an exercice). Moreover, it is known that there are exactly $2^{\ell}-1$ parabolic subgroups $P \supset B_+$, labeled by non-empty subsets of the set of simple roots of the Lie algebra $\mathfrak g$. 

{\bf The diagram of $P$ and of the corresponding parabolic subalgebra $\mathfrak p$},
$$
\mathfrak b_+ \subset \mathfrak p\subset \mathfrak g^{\CC}
$$
is the Dynkin's diagram of the Lie algebra $\mathfrak g$, where the vertices
($=$simple roots) are labeled by $\{0,1\}$. E.g., $P=B_+$, if all simple roots are labeled by $1$.


We have the complex manifold $M=G^{\CC}/P$. The compact Lie group $G$ acts transitively on $M$. We may consider 
$M=G/G\cap P$ as a flag manifold  of $G$ with an $G$-invariant complex structure and the standard basepoint $o=e(G\cap P)$.

Here is an equivalent infinitesimal construction:

\begin{DEF*}
 
A {\bf standard complex form} of a flag manifold $M=G/H$ is a pair $(G/xHx^{-1},\,\mathfrak p)$, where $x\in G$, and $\mathfrak p$ is a parabolic subalgebra of $\mathfrak g^{\CC}$ satisfies the following properties:
$$
\mathfrak b_+ \subset \mathfrak p,\qquad
\mathfrak p\cap\mathfrak g =\opn{Ad}(x)\,\mathfrak h,
\qquad
\mbox{and hence $T^{\ell}\subset xHx^{-1}$.}
$$
The {\bf diagram of a standard complex form} is the diagram of $\mathfrak p$.

\end{DEF*} 
 
Thus, we can associate with a flag manifold $M=G/H$ with a fixed invariant complex structure the unique standard complex form $(G/xHx^{-1},\,\mathfrak p)$ of $M$ and the natural basepoint $xH\in G/H$ (for the uniqueness, cf. the above remark). We obtain:

\begin{CLAIM*}

Any complex flag manifold has a unique standard complex form.
Any flag manifold has at least one standard complex form.

\end{CLAIM*}
 
Assume for simplicity that $(G/H,\mathfrak p)$ is a standard complex form of $M$ (that is, $x\in H$). We describe explicitly the diagram of $(G/H,\mathfrak p)$.
A root $\alpha$ of the Lie algebra $\mathfrak g$ is positive, if $\mathfrak g^\alpha\in \mathfrak b_+ $. In particular let $\gamma$ be every simple root. Then $\gamma$ is labeled
\begin{itemize} 
\item  by $1$,
if $\mathfrak g^\gamma$ belongs to the nilradical of $\mathfrak p $, that is, $\mathfrak g^\gamma\notin\mathfrak h^{\CC}$, 

\item by $0$, otherwise.
\end{itemize}

\smallskip

Given a standard complex form $(M=G/H,\mathfrak p)$ we may easily reconstruct the $T$-root system $\Omega=\Omega_{G/H}$ as follows.
Assume that $R$ is the root system of $\mathfrak g $ with respect to the fixed maximal torus $T^{\ell}$. Let $\{\gamma \}$ be the basis of simple roots, and
$$
R_M = \biggl\{\sum_{\gamma} k_{\gamma} \gamma \in R \,\biggm|\, \sum_{\mbox{$\gamma$ is labeled $1$}} |k_{\gamma}|>0  \biggr\}.
$$
Assume $\alpha_{\gamma}=\gamma+Q_H$, \enskip  cf.~\S\ref{sect:1}. Then
$$
\Omega = \biggl\{\sum_{\mbox{$\gamma$ is labeled $1$}} k_{\gamma} \alpha_\gamma  \,\biggm|\, \sum_{\gamma} k_{\gamma}\gamma \in R_M  \biggr\}.
$$

\smallskip



\smallskip

Thus, we may associate with a standard complex form $(M=G/H,\mathfrak p)$ a basis in $\RR^k$, $k=\opn{rank}(\Omega)$,
namely, the basis $\{\alpha_{\gamma }=\gamma+Q_H : \gamma $ is labeled $1\}$.
This is called {\bf the basis of simple T-roots}. 

The set $\Omega^+$ of {\bf positive T-roots} is the intersection of $\Omega$ with the closed cone in $\RR^k$ generated by simple $T$-roots $\alpha_i$, $i=1,\dots,k$. Obviously, every $\omega\in\Omega$ is either positive ($\omega\in\Omega^+$), or negative ($-\omega\in\Omega^+$).

\medskip

{\bf Remark. } The positive Weyl chamber in $(\RR^k)^* = Lie\,(Z(H))$ associated with a standard complex form of $G/H$ is the simplicial cone
$$
C=\{h \in Lie\,(Z(H)) :  \langle h,\,\, \alpha_i\rangle >0,\, i=1,\dots,k \}.
$$
More generally, Weyl chambers in $(\RR^k)^* = Lie\,(Z(H))$ are the connected components of the open set $\{h:\prod_{\textstyle \omega\in\Omega_{G/H} }\, \langle h,\omega\rangle <0 \}$. Every Weyl chamber is a simplicial cone; there is a one-to-one correspondence between the set of such cones and the set of the invariant complex structures on $G/H$ \cite{AP}. In this context, a standard complex form represents a  class of equivalent invariant complex structures on the flag manifold $G/H$.

\subsection{Properties of $T$-root systems}
Let $\Omega$ be a $T$-root system of rank $k$ in $\RR^k$. A non-empty subset $\Sigma \subset\Omega$ is 
\begin{itemize} 
\item a closed subsystem, if $\alpha,\beta\in \Sigma$, $\alpha+\beta\in\Omega$ implies that $\alpha+\beta\in\Sigma$;
\item a symmetric subsystem, if $\alpha\in \Sigma$ implies that $-\alpha\in\Sigma$;
\item a complete subsystem, if $\Sigma=\Omega\cap V$ for some vector subspace $V\subset \RR^k$.
\end{itemize}
The system $\Omega$ is called 
\begin{itemize} 
\item \textbf{reducible}, if it can be represented by a union $\Omega=\Sigma_1\cup\Sigma_2$ of two complete subsystems $\Sigma_i=\Omega\cap V_i$, $i=1,2$, with $V_1\cap V_2=(0)$,
\item \textbf{irreducible}, otherwise.
\end{itemize}

\begin{LEM}\label{LEM:p-1} 

The group $G$ of a flag manifold $M$ is simple if and only if the $T$-root system $\Omega_M$ is irreducible.

\end{LEM}

\begin{CLAIM*}

Every closed symmetric subsystem $\Sigma\subset\Omega$ is a $T$-root system of rank
$$
\opn{rank}(\Sigma) = \dim(\opn{span}( \Sigma)) 
,\qquad
\opn{rank}(\Sigma)\le \opn{rank}(\Omega),
$$
in the vector space $\opn{span}(\Sigma)$ spanned by $\Sigma$. 

\end{CLAIM*}

Claim follows immediately from the next Lemma~\ref{LEM:p-2} that we state without proof.

Assume that $L$ is a compact connected semisimple Lie group and $N=L/H$ is a flag manifold such that 
$$
\Omega_N=\Omega.
$$

\begin{LEM}\label{LEM:p-2} 

Let $\Omega=\Omega_N$ be $T$-root system of a flag manifold $N=L/H$, and $\Sigma\subset\Omega$ a closed symmetric subsystem. Then there exists a $H$-invariant connected semisimple subgroup $G\subset L$ such that $M=G/G\cap H$ is a flag manifold with the $T$-root system
$$
\Omega_M = \Sigma
$$ 
of rank
$
\opn{rank}(\Omega_M) = \dim(\opn{span}( \Sigma)).
$
Moreover, if $\Sigma$ is irreducible, then the group $G$ of $M=G/G\cap H$ is simple.

\end{LEM}

The last assertion follows from Lemma~\ref{LEM:p-1}.

Note that $G\cap H$ is a normal subgroup of the group $H$. This implies:

\begin{COR*} If $rank(L)-rank(G)=1$, then $[H,H]\subset G\cap H$.
\end{COR*}

\begin{LEM}\label{LEM:p-3} 

Let $\Omega$ be an irreducible system of $T$-roots of rank $k$, 
and $\Sigma_r$ a complete irreducible subsystem of rank $r<k$. Then there exists a complete irreducible subsystem $\Sigma_{r+1}$ of rank $r+1$ such that
$$
\Sigma_{r}\subset\Sigma_{r+1}\subset\Omega.
$$

\end{LEM}

\subsection{Classification of the exceptional flag manifolds}

The classification of exceptional flag manifolds is, of course, very known. It is useful to reading the present paper. Note that ours calculations are \textbf{independent} of this classification. But some checks must be simplified, if we use the following proposition.

\begin{PROP}\label{PROP:f.m.}



Let $M=G/H$ be a flag manifold with $b_2(M)>1$ of an exceptional compact Lie group $G$ of type $E_\ell$, $\ell=6,7,8$. Assume that $H$ contains the fixed maximal torus  $T^\ell$ of $G$. Then up the natural action of the corresponding Weyl group $W=\opn{Norm}_G(T^\ell)/T^\ell$, the group $H$ is of one and only one of the following kinds:
 
\begin{enumerate} 
\item $H$ contains a normal subgroup of type $E_6$ or $D_k$, $k>3$ (there is $1,3,8$ such manifolds $G/H$ respectively),

\item $[H,H]$ belongs to the standard $A_{\ell -1}$ (there is $p(\ell)-1=10,14,21$ such manifolds, where $p(n)$ is the number of partitions of $n$),

\item $H$ can be described by one of the following diagrams:
\\[1ex]
${\sM 0&1&0&1&0\\ &&0\endsM\ig}$ $E_6/T^2\cdot A_1\cdot A_1\cdot A_2$,
\\[1ex]
${\sM 1&0&0&1&0&0\\ &&&0\endsM\ig}$ $E_7/T^2\cdot A_1\cdot A_2\cdot A_2$,\\
${\sM 0&1&0&1&0&0\\ &&&0\endsM\ig}$ $E_7/T^2\cdot A_1\cdot A_1\cdot A_1\cdot A_2$,\\
${\sM 0&0&0&1&1&0\\ &&&0\endsM\ig}$ $E_7/T^2\cdot A_1\cdot A_1\cdot A_3$,\\
${\sM 0&0&0&0&1&1\\ &&&0\endsM\ig}$ $E_7/T^2\cdot [A_5]'' $,\\
${\sM 0&1&0&1&1&0\\ &&&0\endsM\ig}$ $E_7/T^3\cdot A_1\cdot A_1\cdot A_1\cdot A_1 $,\\
${\sM 0&0&0&1&1&1\\ &&&0\endsM\ig}$ $E_7/T^3\cdot [A_1\cdot A_3]'' $,\\
${\sM 0&1&0&1&1&1\\ &&&0\endsM\ig}$ $E_7/T^4\cdot [A_1\cdot A_1\cdot A_1]'' $,
\\[1ex]
${\sM 1&0&0&0&1&0&0\\ &&&&0\endsM\ig}$ $E_8/T^2\cdot A_1\cdot A_2\cdot A_3$,\\
${\sM 1&0&1&0&1&0&0\\ &&&&0\endsM\ig}$ $E_8/T^3\cdot A_1\cdot A_1\cdot A_1\cdot A_2$,\\
${\sM 0&0&1&0&1&0&0\\ &&&&0\endsM\ig}$ $E_8/T^2\cdot A_1\cdot A_1\cdot A_2\cdot A_2$,\\
${\sM 0&1&0&0&0&1&0\\ &&&&0\endsM\ig}$ $E_8/T^2\cdot A_1\cdot A_1\cdot A_4$.

\end{enumerate}

\end{PROP}
Remark that in general a flag manifold of the kind 3) has more than one standard complex forms.


\begin{proof}[Proof for $G=E_8$]
We prove the case 3) using the combinatorics of Dynkin's diagrams. Let $H$ is not of the kinds 1),2). Then it can be described by one of the diagrams:
$$
a){\sM x&x&x&1&x&x&x\\ &&&&0\endsM\ig},
\quad
b){\sM x&x&x&0&0&1&x\\ &&&&0\endsM\ig},
\quad
c){\sM x&x&x&0&1&x&x\\ &&&&0\endsM\ig},
$$
where each $x$ must be replaced by $0$ or $1$. 

a) The condition $a)$ implies
${\sM x&x&x&1&x&0&0\\ &&&&0\endsM\ig}$
(otherwise $H$ reduces to $2)$)
Then
${\sM x&x&0&1&x&0&0\\ &&&&0\endsM\ig}$;
otherwise we reduce $H$ to the case 2)
using $-w\in W(E_8)$, where $-1_8$ is an element of the Weyl group $W(E_8)$,
and $w$ is the element of the maximum length of an appropriate subgroup $W(D_5)$. Hence, we have 
${\sM 0&x&0&1&x&0&0\\ &&&&0\endsM\ig}$. 
Because $b_2(M)>1$, we obtain the diagrams
$$
{\sM 0&0&0&1&1&0&0\\ &&&&0\endsM\ig}
,\quad
{\sM 0&1&0&1&1&0&0\\ &&&&0\endsM\ig}
,\quad
{\sM 0&1&0&1&0&0&0\\ &&&&0\endsM\ig}
$$
of the flag manifolds
$E_8/T^2\cdot A_1\cdot A_2\cdot A_3$,
$E_8/T^3\cdot A_1\cdot A_1\cdot A_1\cdot A_2$,
$E_8/T^2\cdot A_1\cdot A_1\cdot A_4$.

c) The condition $c)$ implies 
${\sM 0&0&x&0&1&x&x\\ &&&&0\endsM\ig}$
(otherwise $H$ reduces to $a)$), and hence
${\sM 0&0&x&0&1&0&0\\ &&&&0\endsM\ig}$
(otherwise $H$ reduces to $a)$ or $b)$).
Since $b_2(M)\ge2$, we obtain the diagram
${\sM 0&0&1&0&1&0&0\\ &&&&0\endsM\ig}$
of $E_8/T^2\cdot A_1\cdot A_1\cdot A_2\cdot A_2$.

b) We prove that $b)$ reduces to $a)$.
The condition $b)$ implies
${\sM 0&x&x&0&0&1&x\\ &&&&0\endsM\ig}$
(otherwise $H$ reduces to 2)), and
${\sM 0&x&0&0&0&1&x\\ &&&&0\endsM\ig}$
(otherwise $H$ reduces to $a)$). Then
${\sM 0&x&0&0&0&1&0\\ &&&&0\endsM\ig}$;
otherwise we reduce $H$ to the case 2)
using $-w\in W(E_8)$, where 
$w$ is the element of the maximum length of an appropriate subgroup $W(D_7)$ (this is correct, since $D_7$ is of an odd rank). Because $b_2(M)>1$, we obtain the diagram
$
{\sM 0&1&0&0&0&1&0\\ &&&&0\endsM\ig}
$
of the flag manifold $E_8/T^2\cdot A_1\cdot A_1\cdot A_4$. Let $u$ be the element of the maximum length of $W(E_6)$. Then
$$
-u:
{\sM 0&1&0&0&0&1&0\\ &&&&0\endsM\ig}
\longrightarrow
{\sM 0&1&0&1&0&0&0\\ &&&&0\endsM\ig},
$$
where the right diagram is of the kind $a)$.
\end{proof}


\begin{COR*} 

There exist $77$ flag manifolds $M$ with $b_2(M)\ge2$ of the exceptional compact simple Lie groups, namely:
\begin{quote} 
$33=8+21+4$ flag manifolds of $E_8$,\\
$24=3+14+7$ flag manifolds of $E_7$,\\
$12=1+10+1$ flag manifolds of $E_6$,\\
$7$ flag manifolds of $F_4$,\\
$1$ flag manifold of $G_2$.
\end{quote}

\end{COR*}

\section{Classification of irreducible T-root systems of rank two}\label{sect:trootwo}

We consider a flag manifold $M=G/H$ as a $G$-manifold (obtained by ignoring the basepoint $eH$).
The isomorpism class $\Omega_M$ of the $T$-root systems $\Omega_{G/xHx^{-1}}$ is well-defined (independent of the choice of basepoints $xH$),
and we call $\Omega_M$ the $T$-root system of the flag manifold $M$, cf.~\S\ref{sect:1}.


Any flag manifold $M$ has finitely many standard complex forms.
Standard complex forms of $M$ are defined in~\S\ref{sect:2.1}.

\begin{THM}\label{THM:trootwo16} 

There exist
\begin{center}
\begin{tabular}{rl}
$16$ & irreducible $T$-root systems $\Omega$ of rank two (up to isomorphism),\\
$30$ & flag manifolds $M$ with $b_2(M)=2$ of the Lie groups $E_6$, $E_7$, $E_8$, $F_4$, and $G_2$,\\
$71$ & standard complex forms of such manifolds.
\end{tabular}
\end{center}
The systems $\Omega$ are pictured in
Table~\ref{table:trootwo-Omega}.
The correspondence between systems $\Omega$, flag manifolds $M$ and standard complex forms of $M$ is given by
Tables~\ref{table:trootwo-flag}
and~\ref{table:trootwo-flag-complex}.

\end{THM}
\begin{THM}\label{THM:trootwo} 

Let $\Omega$ be an irreducible $T$-root system of rank two. Then there exists at least one exceptional compact simple Lie group $L$ other than $G_2$, and a unique flag manifold $M$ of $L$ with a $T$-root system isomorphic to $\Omega$.

\end{THM}

\subsection{Detailed classification}

Fix an irreducible $T$-root system $\Omega$ of rank two. Let 
$\opn{Conv}(\Omega)\subset\RR^2$ be its convex hull, and let $P\subset\RR^2$ be a fundamental parallelogram of the vector lattice generated by the elements of $\Omega$.

\begin{PROP}\label{PROP:da} 

$\Omega$ is uniquely determined by two numbers
$$
d=\frac12|\Omega|,\quad
a=\frac{2\opn{area}(\opn{Conv}(\Omega))}{\opn{area}(P)}.
$$

\end{PROP}

\begin{PROP}\label{PROP:trootwo} 

Let $L$ be a compact Lie group of a type 
$E_6$, $E_7$, $E_8$, or $F_4$, and let $(L/H_i,\,\mathfrak p_i)$, $i=1,2$, are two standard complex flag manifolds of $L$. If the $T$-root systems corresponding to $\mathfrak p_i$
are isomorphic to $\Omega$, then the subgoups $H_i$, $i=1,2$ are conjugate.

\end{PROP}

Hence, there is at most one flag manifold of $L$ with the $T$-root system $\Omega$.

Propositon~\ref{PROP:trootwo} and Theorem~\ref{THM:trootwo} follows from Table~\ref{table:trootwo-flag-complex}.

\input Tables/tabtwo.tex
\clearpage

\subsection{Some inclusions between T-root systems}

Consider consequences of the above classification of the irreducible $T$-root systems of rank two.

\begin{PROP}\label{PROP:i-1} 

Let $\Sigma$ and $\Omega$ are irreducible $T$-root systems,
$\Sigma\subset\Omega$,
$n=\opn{rank}(\Omega)\ge  3$,
$\opn{rank}(\Sigma) =2$,
$d=\frac12|\Sigma|\ge 7$.
Assume $\Sigma\subset\Omega$ is a symmetric closed subsystem which generates a plane in $\RR^n$. Let $M$ and $N$ are flag manifolds of compact Lie groups $G$ and $L$ such that $\Omega_M=\Sigma$ and $\Omega_N=\Omega$. Then 
$$
M=E_7/T^2\cdot H,\quad N=E_8/T^3\cdot H
$$
for one of the following subgroup $H$ of rank $ 5 :$
$$
A_1\cdot A_1\cdot A_1\cdot A_2,\quad   
A_1\cdot A_2\cdot A_2,\quad  
A_1\cdot A_1\cdot A_3,\quad  
A_2\cdot A_3.
$$
The flag manifolds $M$ and $N$ depend only of $\Sigma$. 

\end{PROP}

\begin{proof}[Proof] 
By Lemma~\ref{LEM:p-3} there is an irreducible $T$-root system $\Sigma_3$ of rank $3$ such that $\Sigma \subset \Sigma_3 \subset \Omega$. Further, $G$ is an exceptional simple Lie group other than $G_2$, since $\Sigma$ is irreducible and $d\ge 7$. From Table~\ref{table:trootwo-flag} it follows that either $M=E_7/T^2\cdot H$, or $G=E_8$. By Lemma~\ref{LEM:p-2}, \enskip $G=E_7$, $L=E_8$, $\Omega=\Sigma_3$, and $N=E_8/T^3\cdot H$. Then the flag manifold $N$ has one of the following standard complex forms:
{
\input Tables/inclusion

}%
\noindent
It is easy that the diagrams  in each column defines complex forms of the same flag manifold, and hence $N$ depends only of $\Sigma$.
\end{proof}

Assume now $\Sigma$ is an irreducible $T$-root system of rank two satisfying Proposition~\ref{PROP:i-1}, and $\alpha\in\Sigma$ is a $T$-root such that
$$
3\alpha \in \Sigma.
$$
From Table~\ref{table:trootwo-flag} it follows that $\Sigma $ is the root system
of the flag manifold $E_7/T^2\cdot A_1\cdot A_2\cdot A_2$. Thus, we obtain:

\begin{PROP}\label{PROP:i-2}

$E_8/T^3\cdot A_1\cdot A_2\cdot A_2$  is the unique flag manifold $G/H$ with an irreducible $T$-root system $\Omega$ of rank $\ge 3$
such that $\alpha $ and $3\alpha $ are both $T$-roots for some $\alpha\in\Omega$.

\end{PROP}

\begin{proof}[Proof] 
By Lemma~\ref{LEM:p-3} we have $\alpha,3\alpha\in \Sigma\subset\Omega$ for some complete irreducible subsystem $\Sigma$ with $\opn{rank}(\Sigma)=2$. 
From Table~\ref{table:trootwo-flag} it follows that $d=\frac12|\Sigma|\ge 8$, so $\Sigma$ satisfies Proposition~\ref{PROP:i-1}. The assertion follows.
\end{proof}

\begin{proof}[A direct proof ($G=E_8$)] 
The left part of the following table contains the diagrams of the all mutually non-equivalent invariant complex structures on the flag manifold 
$E_8/T^3\cdot A_1\cdot A_2\cdot A_2$.
The right part of this table contains the all positive roots of the Lie algebra $E_8$ with at least $3$
coefficients $k_{i_j} = 0\pmod3$, $i_j \in [1,8]$, such that $k_{i_1}>0$ (there exist $16$ such roots). 
{\footnotesize
\input Tables/roots16
}%
\noindent
(see, e.g., tables in \cite{Bou}).
This table prove
Proposition~\ref{PROP:i-2} for $G=E_8$.
\end{proof}


\begin{COR*}{} Let $\Omega$ be an irreducible $T$-root system of rank $\ge 3$, and $\omega\in\Omega$ a $T$-root.
If $\frac1k\omega \in \Omega$ then either $k\in\{\pm1,\pm2\}$, or $k\in\{\pm3\}$ and $\Omega$ is the $T$-root system of rank $3$,
corresponding to the flag manifold $E_8/T^3\cdot A_1\cdot A_2\cdot A_2$.
\end{COR*}

\ENDINPUT
\section{}

%% file: Tables/tabtwo.tex




\begin{table*}[htb]{}
\caption{Irreducible T-root systems of rank two}
{
$$
\label{table:trootwo-Omega}
\begin{array}{cccc}
\trooII{1}&\trooII{2}&\trooII{3}&\trooII{4}\\\noalign{\medskip}
\trooII{5}&\trooII{6}&\trooII{7}&\trooII{8}\\\noalign{\medskip}
\trooII{9}&\trooII{10}&\trooII{11}&\trooII{12}\\\noalign{\medskip}
\trooII{13}&\trooII{14}&\trooII{15}&\trooII{16}
\end{array}
$$
}\\[2ex]
\footnotesize
$d=\frac12|\Omega|$; $a:$ the normalized area of the convex hull of $\Omega$ (Proposition~\ref{PROP:trootwo})
\end{table*}



\begin{table*}[htb]{}
\caption{Exceptional flag manifolds with irreducible $T$-root systems of rank two}
\footnotesize
\begin{quote}{}
We enumerate the flag manifolds with $T$-root systems $\Omega$ of rank $2$ of the compact Lie groups of types $E_\ell$, $\ell=6,7,8$, and $F_4$.
Notations:
$n^o_{=}:$ the number of picture (Table~\ref{table:trootwo-Omega});
\\ \quad
$d=\frac12|\Omega|$; 
\\ \quad
$a:$ the normalized area of the convex hull of $\Omega$ (Proposition~\ref{PROP:trootwo});
\\ \quad
$\mathrm{type}(\Omega):$ the type of $\Omega$, e.g., the root system of the type $G_2$ or $BC_2$ (cf. Theorem~\ref{THM:classtypes});
\\ \quad
$m=\max\{k: k \omega \in \Omega,\,\omega \in \Omega\}$;
\\ \quad
$f:$ the number of the standard complex forms of a flag manifold.
\\
We use Dynkin's notations \cite[\S3 (17)]{Dynkin}:\enskip
$A_5'$, $A_5''\subset E_7 $, and $ A_k$, $\widetilde A_k\subset  F_4$, where $k \in \{1,2 \}$. 
\end{quote}
\vskip2ex 
$$
\def\9{}\let\TIL\widetilde
\def\arraystretch{1.4}\label{table:trootwo-flag}
\begin{array}{|r|r|r|l|r|ll|ll|ll|ll|}
\hline
n^o_{=} &d &a & \mathrm{type}(\Omega) &m&{{{}}}
E_6/T^2\bullet&f&E_7/T^2\bullet&f&E_8/T^2\bullet&f&F_4/T^2\bullet&f
\\ \hline
16 &  14 &  48   &       &3&{{{}}}         & &             & &A_1A_1A_2A_2 &2&             &  \\ 
15 &  12 &  40   &       &4&{{{}}}         & &             & &A_1A_2A_3    &4&             &  \\ 
14 &  11 &  34   &       &2&{{{}}}         & &             & &A_1A_1A_4    &4&             &  \\ 
13 &  10 &  32   &       &3&{{{}}}         & &             & &A_2A_4       &4&             &  \\ 
12 &  10 &  28   &       &2&{{{}}}         & &             & &A_3A_3       &2&             &  \\ 
11 &   9 &  28   &       &3&{{{}}}         & &             & &A_1A_5       &3&             &  \\ 
10 &   9 &  24   &       &2&{{{}}}         & &A_1A_1A_1A_2 &1&A_2D_4       &1&A_2          &1 \\ \hline
9  &   8 &  24   &       &2&{{{}}}         & &  A_1A_1A_3  &3&A_1D_5       &3&A_1\TIL A_1  &3 \\ 
8  &   8 &  20   &       &3&{{{}}}         & &  A_1A_2A_2  &3&    \9       & &   \9        &  \\ 
7  &   7 &  20   &       &2&{{{}}}         & &  A_2A_3     &3&A_6          &3&   \9        &  \\ \hline
6  &   6 &  18   &G_2    &1&{{{}}}A_2A_2   &1&  A_5''      &1&E_6          &1&\TIL A_2     &1 \\ 
5  &   6 &  16   &BC_2   &2&{{{}}} \9      & &  A_1D_4     &1&D_6          &1& B_2         &1 \\ 
4  &   6 &  14   &       &2&{{{}}}A_1A_1A_2&5&  A_1A_4     &5&   \9        & &   \9        &  \\ \hline
3  &   5 &  12  &BC_{2,1}&2&{{{}}}A_1A_3   &4&  A_5'       &2&   \9        & &   \9        &  \\ 
2  &   4 &   8   &B_2    &1&{{{}}}A_4      &4&  D_5        &2&   \9        & &   \9        &  \\ 
1  &   3 &   6   &A_2    &1&{{{}}}D_4      &1&  \9         & &   \9        & &   \9        &  \\ \hline
\hline  \multicolumn{5}{|c|}{\text{╒А╔ё╝:}}
                            &5       &15&  9         &21 &  11       &28&  4          &6 \\ \hline
\end{array}
$$
\begin{quote}
$16$ systems of $T$-roots,\enskip $30 = 5+9+11+4+1$ flag manifolds (together with $G_2/T^2$),\enskip $71=15+21+28+6+1$ standard complex forms.
\end{quote}
%
\end{table*}





\begin{table*}[htb]{} 
\caption{Standard complex forms of exceptional flag manifolds with irreducible $T$-root systems of rank two}
\label{table:trootwo-flag-complex}
{\footnotesize
We enumerate the standard complex forms of the flag manifolds with $T$-root systems $\Omega$ of rank $2$ of the compact Lie groups $E_6$, $E_7$, $E_8$, and $F_4$.
}\small
$$
\def\arraystretch{1.4}
\begin{array}{|l|l|l|l|l|l|}
\hline
d&a&E_6&E_7&E_8&F_4
\\ \hline
9&28& & &{\sM 0&1&0&0&0&0&0\\ &&&&1\endsM\ig}&
\\  
  &  & & &{\sM 0&0&0&0&0&1&0\\ &&&&1\endsM\ig}&
\\  
  &  & & &{\sM 1&0&0&0&0&1&0\\ &&&&0\endsM\ig}&
\\  \hline
9&24& &{\sM 0&1&0&1&0&0\\ &&&0\endsM\ig}&{\sM 0&0&1&0&0&0&1\\ &&&&0\endsM\ig}&
{\sM 1&1\Leftarrow0&0\endsM\ig}
\\  \hline
8&24& &{\sM 0&0&0&1&0&1\\ &&&0\endsM\ig}&{\sM 1&0&1&0&0&0&0\\ &&&&0\endsM\ig}&
{\sM 0&1\Leftarrow0&1\endsM\ig}
\\  
  &  & &{\sM 0&0&0&1&1&0\\ &&&0\endsM\ig}&{\sM 0&1&1&0&0&0&0\\ &&&&0\endsM\ig}
&{\sM 0&1\Leftarrow1&0\endsM\ig}
\\  
  &  & &{\sM 0&1&0&0&1&0\\ &&&0\endsM\ig}&{\sM 0&1&0&0&0&0&1\\ &&&&0\endsM\ig}
&{\sM 1&0\Leftarrow1&0\endsM\ig}
\\  \hline
8&20& &{\sM 0&0&1&0&1&0\\ &&&0\endsM\ig}& &
\\  
  &  & &{\sM 0&0&1&1&0&0\\ &&&0\endsM\ig}& &
\\  
  &  & &{\sM 1&0&0&1&0&0\\ &&&0\endsM\ig}& &
\\  \hline
7&20& &{\sM 0&0&0&1&0&0\\ &&&1\endsM\ig}&{\sM 0&0&0&0&0&1&1\\ &&&&0\endsM\ig}&

\\  
  &  & &{\sM 0&0&1&0&0&0\\ &&&1\endsM\ig}&{\sM 0&0&0&0&0&0&1\\ &&&&1\endsM\ig}
&
\\  
  &  & &{\sM 0&0&1&0&0&1\\ &&&0\endsM\ig}&{\sM 1&0&0&0&0&0&0\\ &&&&1\endsM\ig}
&
\\  \hline
6&18&{\sM 0&0&1&0&0\\ &&1\endsM\ig}&{\sM 0&0&0&0&1&1\\ &&&0\endsM\ig}&{\sM 1&
1&0&0&0&0&0\\ &&&&0\endsM\ig}&{\sM 0&0\Leftarrow1&1\endsM\ig}
\\  \hline
6&16& &{\sM 0&1&0&0&0&1\\ &&&0\endsM\ig}&{\sM 1&0&0&0&0&0&1\\ &&&&0\endsM\ig}&
{\sM 1&0\Leftarrow0&1\endsM\ig}
\\  \hline
6&14&{\sM 0&0&1&0&1\\ &&0\endsM\ig}&{\sM 0&1&0&0&0&0\\ &&&1\endsM\ig}& &
\\  
  &  &{\sM 0&0&1&1&0\\ &&0\endsM\ig}&{\sM 0&0&0&0&1&0\\ &&&1\endsM\ig}& &
\\  
  &  &{\sM 0&1&0&1&0\\ &&0\endsM\ig}&{\sM 1&0&0&0&1&0\\ &&&0\endsM\ig}& &
\\  
  &  &{\sM 0&1&1&0&0\\ &&0\endsM\ig}&{\sM 1&0&1&0&0&0\\ &&&0\endsM\ig}& &
\\  
  &  &{\sM 1&0&1&0&0\\ &&0\endsM\ig}&{\sM 0&1&1&0&0&0\\ &&&0\endsM\ig}& &
\\  \hline
%
%
5&12&{{\sM 0&0&0&1&0\\ &&1\endsM\ig}\quad  {\sM 0&1&0&0&0\\ &&1\endsM\ig}}&{\sM 0&0&
0&0&0&1\\ &&&1\endsM\ig}& &
\\  
  &  &{{\sM 1&0&0&1&0\\ &&0\endsM\ig} \quad {\sM 0&1&0&0&1\\ &&0\endsM\ig}}&{\sM 1&0&
0&0&0&0\\ &&&1\endsM\ig}& &
\\  \hline
4&8&{{\sM 1&1&0&0&0\\ &&0\endsM\ig} \quad  {\sM 0&0&0&1&1\\ &&0\endsM\ig}}&{\sM 1&0&0&
0&0&1\\ &&&0\endsM\ig}& &
\\  
  &  &{{\sM 1&0&0&0&0\\ &&1\endsM\ig}\quad {\sM 0&0&0&0&1\\ &&1\endsM\ig}}&{\sM 1&1&
0&0&0&0\\ &&&0\endsM\ig}& &
\\  \hline
3&6&{\sM 1&0&0&0&1\\ &&0\endsM\ig}& & &
\\  \hline
\end{array}
\qquad \qquad
\def\arraystretch{1.4}
\begin{array}{|l|l|l|}
\hline
d&a&E_8
\\ \hline
14&48&{\sM 0&0&1&0&1&0&0\\ &&&&0\endsM\ig}
\\  
  &  &{\sM 0&1&0&0&1&0&0\\ &&&&0\endsM\ig}
\\  \hline
12&40&{\sM 0&0&1&0&0&1&0\\ &&&&0\endsM\ig}
\\  
  &  &{\sM 0&0&0&1&0&1&0\\ &&&&0\endsM\ig}
\\  
  &  &{\sM 0&0&0&1&1&0&0\\ &&&&0\endsM\ig}
\\  
  &  &{\sM 1&0&0&0&1&0&0\\ &&&&0\endsM\ig}
\\  \hline
11&34&{\sM 0&0&0&0&1&0&1\\ &&&&0\endsM\ig}
\\  
  &  &{\sM 0&0&0&0&1&1&0\\ &&&&0\endsM\ig}
\\  
  &  &{\sM 0&1&0&0&0&1&0\\ &&&&0\endsM\ig}
\\  
  &  &{\sM 0&1&0&1&0&0&0\\ &&&&0\endsM\ig}
\\  \hline
10&32&{\sM 0&0&0&0&1&0&0\\ &&&&1\endsM\ig}
\\  
  &  &{\sM 0&0&1&0&0&0&0\\ &&&&1\endsM\ig}
\\  
  &  &{\sM 0&0&1&1&0&0&0\\ &&&&0\endsM\ig}
\\  
  &  &{\sM 1&0&0&1&0&0&0\\ &&&&0\endsM\ig}
\\  \hline
10&28&{\sM 0&0&0&1&0&0&0\\ &&&&1\endsM\ig}
\\  
  &  &{\sM 0&0&0&1&0&0&1\\ &&&&0\endsM\ig}
\\  \hline
\end{array}
$$
\end{table*}



%% file: Tables/inclusion.tex

$$
\begin{array}{cccc}
A_1A_1A_1A_2 &  A_1A_2A_2 & A_1A_1A_3 & A_2A_3 
\\[1ex]
{
\def\arraystretch{1.4}
\begin{array}{|l|}
\hline
{\sM 1&0&1&0&1&0&0\\ &&&&0\endsM\ig}
\\  \hline
{\sM 0&1&1&0&1&0&0\\ &&&&0\endsM\ig}
\\  \hline
{\sM 0&1&0&1&1&0&0\\ &&&&0\endsM\ig}
\\  \hline
{\sM 0&1&0&1&0&1&0\\ &&&&0\endsM\ig}
\\  \hline
{\sM 0&1&0&0&1&1&0\\ &&&&0\endsM\ig}
\\  \hline
{\sM 0&1&0&0&1&0&1\\ &&&&0\endsM\ig}
\\  \hline
{\sM 0&0&1&0&1&1&0\\ &&&&0\endsM\ig}
\\  \hline
{\sM 0&0&1&0&1&0&1\\ &&&&0\endsM\ig}
\\  \hline
\end{array}
}&
{
\def\arraystretch{1.4}
\begin{array}{|l|}
\hline
{\sM 1&1&0&0&1&0&0\\ &&&&0\endsM\ig}
\\  \hline
{\sM 1&0&0&1&1&0&0\\ &&&&0\endsM\ig}
\\  \hline
{\sM 1&0&0&1&0&1&0\\ &&&&0\endsM\ig}
\\  \hline
{\sM 0&1&0&0&1&0&0\\ &&&&1\endsM\ig}
\\  \hline
{\sM 0&0&1&1&1&0&0\\ &&&&0\endsM\ig}
\\  \hline
{\sM 0&0&1&1&0&1&0\\ &&&&0\endsM\ig}
\\  \hline
{\sM 0&0&1&0&1&0&0\\ &&&&1\endsM\ig}
\\  \hline
{\sM 0&0&1&0&0&1&0\\ &&&&1\endsM\ig}
\\  \hline
\end{array}
}&
{
\def\arraystretch{1.4}
\begin{array}{|l|}
\hline
{\sM 1&0&1&0&0&1&0\\ &&&&0\endsM\ig}
\\  \hline
{\sM 1&0&0&0&1&1&0\\ &&&&0\endsM\ig}
\\  \hline
{\sM 1&0&0&0&1&0&1\\ &&&&0\endsM\ig}
\\  \hline
{\sM 0&1&1&0&0&1&0\\ &&&&0\endsM\ig}
\\  \hline
{\sM 0&1&0&1&0&0&1\\ &&&&0\endsM\ig}
\\  \hline
{\sM 0&1&0&1&0&0&0\\ &&&&1\endsM\ig}
\\  \hline
{\sM 0&1&0&0&0&1&0\\ &&&&1\endsM\ig}
\\  \hline
{\sM 0&0&0&1&1&1&0\\ &&&&0\endsM\ig}
\\  \hline
{\sM 0&0&0&1&1&0&1\\ &&&&0\endsM\ig}
\\  \hline
{\sM 0&0&0&1&0&1&0\\ &&&&1\endsM\ig}
\\  \hline
\end{array}
}&
{
\def\arraystretch{1.4}
\begin{array}{|l|}
\hline
{\sM 1&0&0&1&0&0&1\\ &&&&0\endsM\ig}
\\  \hline
{\sM 1&0&0&1&0&0&0\\ &&&&1\endsM\ig}
\\  \hline
{\sM 1&0&0&0&1&0&0\\ &&&&1\endsM\ig}
\\  \hline
{\sM 0&0&1&1&0&0&1\\ &&&&0\endsM\ig}
\\  \hline
{\sM 0&0&1&1&0&0&0\\ &&&&1\endsM\ig}
\\  \hline
{\sM 0&0&1&0&0&1&1\\ &&&&0\endsM\ig}
\\  \hline
{\sM 0&0&1&0&0&0&1\\ &&&&1\endsM\ig}
\\  \hline
{\sM 0&0&0&1&1&0&0\\ &&&&1\endsM\ig}
\\  \hline
{\sM 0&0&0&1&0&1&1\\ &&&&0\endsM\ig}
\\  \hline
{\sM 0&0&0&1&0&0&1\\ &&&&1\endsM\ig}
\\  \hline
\end{array}
}
\end{array}
$$

%% file: Tables/roots16.tex
$$
\let\meti\null
\def\arraystretch{1.4}
\begin{array}{|c||cc|}
\hline
\mbox{diagram}&\mbox{root 1}&\mbox{root 2}\\ \hline
\meti
	{\sM 1&1&0&0&1&0&0\\ &&&&0\endsM\ig}
&
	{\sM 0&0&1&2&3&2&1\\ &&&&1\endsM\ig}&
\meti
	{\sM 0&0&1&2&3&2&1\\ &&&&2\endsM\ig}
\\ \hline
\meti
	{\sM 1&0&0&1&1&0&0\\ &&&&0\endsM\ig}
&
	{\sM 0&1&2&3&3&2&1\\ &&&&1\endsM\ig}&
\meti
	{\sM 0&1&2&3&3&2&1\\ &&&&2\endsM\ig}
\\  \hline
\meti
	{\sM 1&0&0&1&0&1&0\\ &&&&0\endsM\ig}
&
	{\sM 0&1&2&3&4&3&1\\ &&&&2\endsM\ig}&
	{\sM 0&1&2&3&4&3&2\\ &&&&2\endsM\ig}
\\  \hline
\meti
	{\sM 0&0&1&1&1&0&0\\ &&&&0\endsM\ig}
&
	{\sM 1&2&3&3&3&2&1\\ &&&&1\endsM\ig}&
\meti
	{\sM 1&2&3&3&3&2&1\\ &&&&2\endsM\ig}
\\  \hline
\meti
	{\sM 0&0&1&1&0&1&0\\ &&&&0\endsM\ig}
&
	{\sM 1&2&3&3&4&3&1\\ &&&&2\endsM\ig}&
\meti
	{\sM 1&2&3&3&4&3&2\\ &&&&2\endsM\ig}
\\  \hline
\meti
	{\sM 0&0&1&0&0&1&0\\ &&&&1\endsM\ig}
&
	{\sM 1&2&3&4&5&3&1\\ &&&&3\endsM\ig}&
	{\sM 1&2&3&4&5&3&2\\ &&&&3\endsM\ig}
\\  \hline
\meti
	{\sM 0&0&1&0&1&0&0\\ &&&&1\endsM\ig}
&
	{\sM 1&2&3&4&6&4&2\\ &&&&3\endsM\ig}&
	{\sM 1&2&3&5&6&4&2\\ &&&&3\endsM\ig}
\\  \hline
\meti
	{\sM 0&1&0&0&1&0&0\\ &&&&1\endsM\ig}
&
	{\sM 1&3&4&5&6&4&2\\ &&&&3\endsM\ig}&
	{\sM 2&3&4&5&6&4&2\\ &&&&3\endsM\ig}
\\  \hline
\end{array}
$$

%% file: T-1.tex
\input hat
\NNULL{\setcounter{section}{3} 
\theoremstyle{plain}
\newtheorem*{CLAIM*}{Claim}
}

\section{$T$-root systems of exceptional flag manifolds with $b_2\ge2$}\label{sect:4}

In this paper we classify the $T$-root systems of ranks $\ge 2$ corresponding to flag manifolds of exceptional simple Lie groups. 
We state main results of the classification.

\begin{THM}\label{THM:L/H_i} Let $L$ be an exceptional compact simple Lie group,
$L/H_i$, $i=1,2$ are flag manifolds with second Betti numbers $b_2(L/H_i)=k$, and $\Omega_i$ the corresponding $T$-root systems of rank $k$.  If $k\ge2$, then the following conditions are equivalent:
\begin{enumerate} 
\item $\Omega_1$ and $\Omega_2$ are isomorphic;
\item $H_1$ and $H_2$ are conjugate in $L$.
\end{enumerate}
\end{THM}

\begin{proof}[Sketch of a Proof] Let $\varSigma$ be the set of Dynkin diagrams of standard complex forms of the flag manifolds $M$ of the group $L$ with $b_2>1$.
Each vertex of $S\in\varSigma$ is labelled $0$ or $1$. 
 
1) Given every diagram $S\in\varSigma$, we calculate the corresponding $T$-root system $\Omega=\Omega(S)$ (as in~\S\,\ref{sect:2.1}), and four numeric invariants $k,d,c,v$ of $\Omega$, where $k=\opn{rank}(\Omega)$, $d=\frac12|\Omega|$, and $c,v$ are defined in \S\,\ref{sect:tables} below. 

2) One can check directly that the diagrams $S$ with fixed $k,d,c,v$ correspond to conjugated flag manifolds of $L$ (under the group of inner automorphisms of $L$). 

Remark that 1) can be easily programmed.

\footnotesize (Remark 2.
For another proof we may use Proposition~\ref{PROP:f.m.} and associate with each $M$ a standard complex form $S_M\in\varSigma$. Then it is sufficient to check that $S_M$ is uniquely determined by $k,d,c,v$.)
\end{proof}

\begin{THM}\label{THM:EF} 
The classification of systems $\Omega$ of $T$-roots
of the flag manifolds $M$ with 
$b_2(M)\equiv\opn{rank}(\Omega)\ge2$ of groups $E_6$, $E_7$, $E_8$, and $F_4$ is given by Table~\ref{table:31} in \S\,\ref{sect:tables} below.
We skip the group $G_2$. The $T$-roots of $M=G_2/T^2$ coincides, obviously, with roots of $G_2$. 
 


\end{THM}

Each line of tables corresponds to a $T$-root system (up to isomorpism).

Now, we give some consequences of such classification. 

\subsection{Flag manifolds with isomorphic $T$-root systems}

By definition of a flag manifold $M=G/H$ we ignore the basepoint $eH$, and consider $M$ with the basepoints $xH$. Recall that the $T$-root system $\Omega$ of a coset space $G/xHx^{-1}$ is unique, up to isomorphism  (cf.~\S\ref{sect:1}).

\begin{DEF*}
Given an irreducible $T$-root system $\Omega$, denote by $\mathcal F(\Omega)$
the set of the flag manifolds with $T$-root systems isomorphic to $\Omega$.
\end{DEF*}
\smallskip

\begin{THM}\label{THM:FOm} 

Let $\mathcal F(\Omega)$ contains an exceptional flag manifold, $|\mathcal F(\Omega)|\ge2$, and $\opn{rank}(\Omega)\ge3$. Then $\mathcal F(\Omega)$ is one of the followings\footnote{In Theorem~\ref{THM:FOm} we use Dynkin's notation $[\,]''$, $[\,]'$ for regular semisimple subgroups (subalgebras) of $E_\ell$, $\ell\in\{7,8\}$. Namely, we write $[H]'$ if $H\subset A_\ell$, and $[H]''$ otherwise. In 2), one can write $[4A_1]''$ without $4A_1$, since $4A_1\supset[3A_1]''$; but $E_7/T^3[4A_1]'$ is a non-flag manifold.}:
{\small
\begin{align*}{}
&1)\,&&E_7/T^4[A_1A_1A_1]'',\,&&  E_8/T^4D_4,\,&&  F_4/T^4,&&rank(\Omega _1 )=4,&&d=24;\\
&2)\,&&E_7/T^3A_1A_1A_1A_1,\,&&  E_8/T^3A_1D_4,\,&&  F_4/T^3A_1,&&rank(\Omega _2 )=3,&&d=16;\\
&3)\,&&E_7/T^3[A_1A_3]'',\,&&  E_8/T^3D_5,\,&&  F_4/T^3\widetilde A_1,&&rank(\Omega _3 )=3,&&d=13;\\
&4)\,&&E_6/T^3A_1A_2,\,&& E_7/T^3A_4,&&&&rank(\Omega _4 )=3, &&d=10;\\
&5)\,&&C_{n}/T^3A_{n_1-1}A_{n_2-1}A_{n_3-1},\,&&D_{n+3}/T^3A_{n_1}A_{n_2}A_{n_3},\,&&E_7/T^3D_4,\,&&rank(\Omega _5 )=3,&&d=9;\\
&6)\,&&&&D_{n+3}/T^3A_{n_1}A_{n_2},\, &&E_6/T^3A_3,\,&&rank(\Omega _6 )=3,&&d=8.
\end{align*}
}%
(where $n_i\ge1$, $\sum n_i = n$).
Each of 1)-6) is one of $\mathcal F(\Omega )$.
 

 \end{THM}

Hence, $\Omega_1$ is the exceptional root system of the type $F_4$, since $F_4/T^4\in \mathcal F(\Omega_1 )$.

\begin{proof}[Proof]
The lines 1)-4) are contained in Tables~\ref{table:31} and~\ref{table:Numberings} below, see~\S\,\ref{sect:tables}. Then there exist four sets $\mathcal F(\Omega)$ such that every of them contains at least two exceptional flag manifolds.
The lines 5) and 6) follows from Proposition~\ref{PROP:C_3} below. To complete the proof, we use Proposition~\ref{PROP:excep}. 
\end{proof}

\begin{COR*} 

The condition of Theorem~\ref{THM:FOm} holds for $15$ irreducible $T$-root system $\Omega$ of ranks $\ge2$. 
In particular, $|\mathcal F(\Omega)|\ge2$ for the $T$-root system $\Omega$ of each flag manifold $E_\ell/T^k\cdot H$, where $k\ge2$, and $H$ has a normal subgroup $D_4$, $D_5$, $D_6$, or $E_6$. 

\end{COR*}

We describe now the $T$-root systems 4)-6) of Theorem~\ref{THM:FOm}.

\begin{PROP}{}\label{PROP:C_3} 

Flag manifolds $E_6/T^3\cdot A_3$ and $E_7/T^3\cdot D_4$ have $T$-root systems of types $C_{3,2}$ and $C_3$ respectively.

\end{PROP}

\begin{PROP}\label{PROP:B_3+}

The common  $T$-root system $\Omega$ of $E_6/T^3\cdot A_1\cdot A_2$ and $E_7/T^3\cdot A_4$ is the union of classical root system of type $B_3$ and two opposite vectors $\mathbf v$, $-\mathbf v$, where $\mathbf v$ is the sum of three short positive roots of $B_3$.

\end{PROP}

\begin{proof}[Proof] 
To find $T$-roots, we use each of two diagrams 
${\sM 1&2&3&0&0\\ &&0\endsM\ig}$ and
${\sM 1&2&3&0&0&0\\ &&&0\endsM\ig}$,
where we indicate simple $T$-roots $\alpha_i$, $i=1,2,3$. The major $T$-root  with respect to the basis $\{\alpha_i\}$ is obviously $\mathbf v=\alpha_1+2\alpha_2+3\alpha_3$. It is easy to check that $\Delta=\Omega \setminus \{\mathbf v,-\mathbf v\}$ is the root system of the type $B_3$, and $\alpha_i$, $i=1,2,3$ form the standard basis of $\Delta$ with the diagram 
${\sM 1&2&\Rightarrow&3\endsM\ig}$.
\end{proof}

\begin{EXAM*} 
 
Let $\Omega$ be the $T$-root system of the flag manifold $M=E_7/T^3\cdot D_4$. From Table~\ref{table:31} (\S\,\ref{sect:tables}) it follows that $\Omega$ consists of $2d=18$ vectors in $\RR^3$, and the convex hull of $\Omega$ has $6$ vertices. Obviously, we may write these vertices as $\pm\alpha_1,\pm\alpha_2,\pm\alpha_3$, so $\opn{Conv}(\Omega)$ is an octahedron. The others $T$-roots $\gamma_i$, $i=1,\dots,12$, are the midpoints of twelve edges. 
Therefore $\Omega$ is the root system of the type $C_3$. 
Let $V_2$ be a $2$-subspace of $\RR^3$ generated by midpoints of two intersecting edges, for definiteness, $\frac12(\alpha_1-\alpha_2)$ and $\frac12(\alpha_3-\alpha_2)$. Then $\Omega\cap V_2$ be a hexagonal root system of the type $A_2$. There is the root system of a flag submanifold $M_2=E_6/T^2\cdot \mathrm{Spin}(8)\subset M$, since
$$
T^2\cdot \mathrm{Spin}(8) \subset E_6\subset E_7.
$$
We describe the complexified isotropy representations of $M$ and $M_2$.
The $T^3\cdot \mathrm{Spin}(8)$ module
$\CC\otimes (\mathfrak{e}_7/\RR^3+\mathfrak{spin}(8)$ decomposes as a direct sum of $18$ mutually non-equivalent irreducible submodules $\mathfrak m_{\omega}$, $\omega\in \Omega$.
Since $M$ has a unique standard complex form (with the diagram ${\sM 1&1&0&0&0&1\\ &&&0\endsM\ig}$), the Weyl group $W(M)$ (i.e., $\mathrm{Norm}(H)/H$, $H=T^3\cdot D_4$) acts transitively on the set of Weyl chambers of $\Omega$, and $W(M)=W$, the Weyl group of the root system $C_3$.
So $W(M)$ acts transitively on the subset of modules $\mathfrak m_{\gamma}$
corresponding to midpoints of edges $\gamma=\gamma_i$, $i=1,\dots,12$.
Hence, $\dim_{\CC}(\mathfrak{m}_{\gamma})=1/6\dim_{\RR}(M_2)=8$.
Similarly, $W(M)$ acts transitively on the subset of  the modules $\mathfrak m_{\alpha}$ corresponding to the vertices $\alpha\in\{\pm\alpha_i,i=1,2,3\}$ of octahedron, and hence $\dim_{\CC}(\mathfrak m_{\alpha})=1$, since
$$
\dim_{\RR}(M)-12\cdot 8 = (133 - 3 - 28)-96 = 6.
$$
We may regard the \textbf{real} isotropy module $\mathfrak{e}_6/\RR^2+\mathfrak{spin}(8)$ of $M_2$ as the direct sum of the complex vector and two complex semispinor modules of $\mathrm{Spin}(8)$ and the action of the subgroup $S_3=W_{A_2}=W(M_2)$ 
(permuting the $\alpha_i$, $i=1,2,3$)
of the group $W$ as the Cartan's triality. A proof use symmetries of Dynkin's diagrams $E_6$ and $D_4$. Another proof follows from the diagram
$$
\begin{CD}
\mathrm{SO}(10)/\mathrm{SO}(2)\times\mathrm{SO}(8) @>>>
E_6/T^2\cdot\mathrm{Spin}(8) @>>> 
E_6/T^1\cdot\mathrm{Spin}(10) @.=@.
(\CC\otimes\mathbb C\mathrm a)P^2.
\\
{\sM 1&0&0&0\\ &&0\endsM\ig} @.
{\sM 1&0&0&0&1\\ &&0\endsM\ig} @.
{\sM 0&0&0&0&1\\ &&0\endsM\ig}
\end{CD}
$$
Thus, all $\mathfrak m_\alpha$ are trivial $\mathrm{Spin}(8)$-modules, and 
(given $i<j\le3$) each 
$\mathfrak m_{\gamma}$ is equivalent to $\mathfrak m_{\frac12(\alpha_i-\alpha_j)}$  as a $\mathrm{Spin}(8)$-module
if and only if $\gamma=\frac12(\pm\alpha_i\pm\alpha_j)$.
%
%

\end{EXAM*}

\smallskip

{\bf Note. } We give diagrams of the all standard complex forms of the spaces 1) and 3):
$$
{\sM 0&1&0&1&1&1\\ &&&0\endsM\ig},\quad
{\sM 1&1&1&0&0&0&1\\ &&&&0\endsM\ig},\quad
{\sM 1&1\Leftarrow1&1\endsM\ig}; 
$$
$$
{\sM 0&1&0&0&1&1\\ &&&0\endsM\ig},\quad {\sM 0&0&0&1&1&1\\ &&&0\endsM\ig},\quad
{\sM 1&1&0&0&0&0&1\\ &&&&0\endsM\ig},\quad {\sM 1&1&1&0&0&0&0\\ &&&&0\endsM\ig},\quad
{\sM 1&0\Leftarrow1&1\endsM\ig},\quad {\sM 0&1\Leftarrow1&1\endsM\ig}.
$$
This is a delicate case where exist
two pairs of  distinct flag manifolds of the same form $G/H$ with $11$ standard complex form in each pair. 
Namely,
\begin{center}  
$E_7/T^4[A_1A_1A_1]''$ ($d=24$) and $E_7/T^4[A_1A_1A_1]'$ ($d=23$)
\end{center} 
have $11=1+10$ standard complex forms. Further,
\begin{center} 
$E_7/T^3[A_1A_3]''$ ($d=13$) and $E_7/T^3[A_1A_3]'$ ($d=12$) 
\end{center}
have also $11=2+9$ standard complex forms. The second flag manifold in each pair has proportional $T$-roots $\alpha$ and $2\alpha$.
Remark also that 
\begin{center}
$E_7/T^3[A_1A_3]''$ and $E_7/T^3A_2A_2$
\end{center}
have both $T$-root systems $\Omega$ with $d=\frac12|\Omega|=13$. But the number of vertices of $\opn{Conv}(\Omega)$ equals $12$ and $18$ respectively. 
 
\smallskip

\subsection{Description of $T$-root systems of ranks $k\ge2$ with $d\le10$}

\begin{THM}\label{THM:d<11} 
There are $16$ irreducible $T$-root systems $\Omega$ with $\opn{rank}(\Omega)\ge 2$ and $d=\frac12|\Omega|\le10$,
corresponding to flag manifolds of exceptional Lie groups:
\begin{itemize} 
\item some of $T$-root system of rank two, pictured in Table~\ref{table:trootwo-Omega};
\item the root system of the type $C_3$,
and the $T$-root system of the type $C_{3,2}$, that is, $C_3$ without two long opposite roots $\alpha$ and $-\alpha$ (Proposition~\ref{PROP:C_3});
\item the system $B_3\cup\{\mathbf v,-\mathbf v\}$ of rank $3$, described in Proposition~\ref{PROP:B_3+}.
\end{itemize}

\end{THM}

\subsection{Exceptional and non-exceptional systems of $T$-roots}

\begin{DEF*} Let $\Omega$ be an irreducible $T$-root system. We call $\Omega$ \textbf{non-exceptional} if it corresponds to some flag manifold of a classical compact Lie group, and \textbf{exceptional}, otherwise.
\end{DEF*}

E.g., flag manifolds $E_6/T^3\cdot A_3$ and $E_7/T^3\cdot D_4$ have non-exceptional $T$-root systems of types $C_{3,2}$ and $C_3$ respectively by Proposition~\ref{PROP:C_3}.


\smallskip

From Table~\ref{table:31} in~\S\,\ref{sect:tables} we obtain\,:

\begin{PROP}{}\label{PROP:excep} Let $\Omega=\Omega_k $ be the system of $T$-roots of rank $k$ corresponding to a flag manifold $M$ with $b_2(M)=k$ of a compact Lie group  $L=E_6$, $E_7$, $E_8$, or $F_4$.
Let $M\ne E_6/T^3\cdot A_3$, $M\ne E_7/T^3\cdot D_4$, and $k\ge 3$. Then $\Omega_k$ is exceptional.
\end{PROP}

\begin{proof}[Proof] We may assume that $k\in[3,7]$.  Let $d=\frac12|\Omega_k|$, $i=\frac12|\omega\in\Omega_k:2\omega\in\Omega_k|$.  Assume  $2v\le 2d$ is the number of vertices of the convex hull of $\Omega_k$. From tables we obtain that there one of the following conditions holds.
\begin{enumerate}
\item $k=3$, $d\in[k^2-1,k^2]=[8,9]$, $v<d$; then $M= E_6/T^3\cdot A_3$ or $E_7/T^3\cdot D_4$; 
\item $k\in[4,6]$, $d\in[k^2-1,k^2]$, $v=d$ (there is a unique such system $\Omega_k$ for each $k\in[4,6]$, in particular, the root system $E_6$, $k=6$);
\item $k\in[3,5]$, $d\in[k^2+1,k^2+k]$, $i\in[0,1]$, $i<d-k^2$ (there are $6-k$ such systems $\Omega_k$); 
\item $d>k^2+k$; then $\Omega_k$ is necessarily exceptional.
\end{enumerate}
Consider the cases 2) and 3).
Let $\Omega_k$ is non-exceptional.
 
2)~If $d=k^2$, then $\Omega_k$ is the root system of the type $B_k$ or $C_k$, and vertices of $\opn{Conv}(\Omega_k)$ are the long roots; hence $v<d$. 
If $d=k^2-1$, then $\Omega_k$  is of the type $C_{k,k-1}$, so $\Omega_k$ contains the subsystem $C_{k-1}$; hence $v<d$. But $v=d$. The contradiction.  

3)~If $d\in[k^2+1,k^2+k]$, then $\Omega_k$ is the root system of the type $BC_{k,i}$ (where $BC_{k,k}=BC_k$), and $i=d-k^2$. But $i<d-k^2$. The contradiction.
This prove Proposition~\ref{PROP:excep} and completes the proof of Theorem~\ref{THM:FOm}.
\end{proof}

{\bf Remark. }
The complete list of non-exceptional irreducible $T$-root systems corresponding to flag manifolds of exceptional compact Lie groups is 
$$
A_1, BC_1, A_2, B_2, BC_{2,1}, BC_2, C_{3,2}, C_3,
$$
as a consequence of Propositions~\ref{PROP:excep} and~\ref{PROP:C_3}.

\section{Tables of T-root systems and the corresponding flag manifolds of exceptional simple Lie groups}\label{sect:tables}
 
\subsection{Some invariants of T-root system} 
Fix an irreducible $T$-root system $\Omega$, and the set $\Omega^+$ of the positive $T$-roots. 
Assume  
$$
\Omega ^+_V:= \{\omega \,\in\, \Omega ^+ :
2\omega - \gamma  \notin \Omega^+, \,\,\forall\,\gamma \in \Omega \cup (0)
,\gamma \ne \omega\}. 
$$
For every $\omega\in\Omega$ we get
$g_{\omega }=\max\{k\in\ZZ: \tfrac1k\,\omega \in \Omega \}$, and consider the following numbers:
$$
\begin{aligned}
d&=\frac12\,|\Omega |=|\Omega ^+|,&
v&=|\Omega ^+_V|,\\[1ex]
c&=\prod_{\omega \, \in \,\Omega ^+} g_{\omega },&
w&=\prod_{\omega \, \in\, \Omega ^+_V} g_{\omega }.
\end{aligned}
$$
It is clear that the numbers $c,d,v,w$ are invariant under the action of the group $\opn{Aut}(\Omega)$.

\smallskip

Note that $2v$ is the number of the vertices of $\opn{Conv}(\Omega)$ by the following proposition:

\begin{PROP}\label{PROP:vertices}
The set of the vertices of the convex hull of $\Omega$ coincides with
$$
\Omega_V:= \{\omega \in \Omega :
2\omega - \gamma  \notin \Omega \,\,\forall\,\gamma \in \Omega \cup (0)
,\gamma \ne \omega\},
$$
i.e., with $\{+\omega, -\omega : \omega\in\Omega^+_V \}$. 
 \end{PROP}

Proposition~\ref{PROP:vertices} is not used in the classification of $T$-root systems. 
 
\medskip


We may clarify the meaning of $c$, $w$ using the classification of the $T$-root systems of rank $k=2$.
%
Let 
$$
i(n)=|\omega\in\Omega^+:n\omega\in\Omega^+|,\qquad
j(n)=|\omega\in\Omega^+:n\omega\in\Omega^+_V|.
$$ 
 

\begin{PROP}\label{PROP:ij-2}

Let $\Omega $ be an irreducible $T$-root system with rank $k\ge 2$
of a flag manifold $M\ne E_8/T^2\cdot A_1\cdot A_2\cdot A_3$.
Then $g_{\omega } \in \{1,2,3\}$ for all $\omega \in \Omega $, and
$$ 
c=2^{i(2)}3^{i(3)},\quad
w=2^{j(2)}3^{j(3)}.
$$
\end{PROP}


\begin{PROP}\label{PROP:ij-3}

Let $\Omega $ be an irreducible $T$-root system of a rank $k\ge 3$
of a flag manifolg $M\ne E_8/T^3\cdot A_1\cdot A_2\cdot A_2$.
Then $g_{\omega } \in \{1,2\}$ for any $\omega \in \Omega $, and
$$ 
c=2^{i(2)},\quad
w=2^{j(2)}.
$$
 
\end{PROP}

\begin{proof}[Proof of Propositions~\ref{PROP:ij-2} and~\ref{PROP:ij-3}] 
Case $k=2$. It follows from the classification of $T$-root systems of rank two, that $m:=\max g_{\omega}\le 3$. Cf. Tables~\ref{table:trootwo-flag} and~\ref{table:trootwo-Omega}. 
Case $k\ge3$. Then $m\le2$ by  Proposition~\ref{PROP:i-2} (Corollary).
\end{proof}

\smallskip

To calculate $c$ and $w$ we use the following fact:

\begin{CLAIM*}\label{LEM:gcd} 
$g_{\omega}$ is the g.c.d. of coefficients of $\omega\in\Omega$ with respect to the basis of simple $T$-roots.
\end{CLAIM*}

\subsection{Classification results. Tables}
In the following tables any  $T$-root system $\Omega$ is replaced by its invariants $k,d,c,v,w$, where $k=\opn{rank}(\Omega)$. 

Also as a result of classification of systems $\Omega$ we state the following theorem. 

Let $M_i$, $i=1,2$ are flag manifolds with $b_2(M)\ge 2$ of exceptional compact simple Lie groups $L_i$, $i=1,2$, and $\Omega_i=\Omega_{M_i}$ the corresponding systems of $T$-roots.

\begin{THM}\label{THM:kcdvw} 


The $T$-root systems $\Omega_i=\Omega_{M_i}$, $i=1,2$ with the same invariants $k,d,c,v$ (or the same $k,d,v,w$) are isomorphic.

\end{THM}

\begin{proof}[Proof] The proof is similar to the proof of Theorems~\ref{THM:L/H_i} and~\ref{THM:EF}.
\end{proof} 


Propositions~\ref{PROP:vertices}, \ref{PROP:ij-2} and~\ref{PROP:ij-3} are not used in the classification of $T$-root systems and the proof of Theorem~\ref{THM:kcdvw}.
 
\begin{COR*} 

Theorem~\ref{THM:kcdvw} remains valid for flag manifolds $M_i$ with $b_2(M_i)\ge 1$ of simple compact Lie groups $L_i$, $i=1,2$ (even classical).

\end{COR*}

\begin{proof}[Proof] The proof is similar to the proof of Proposition~\ref{PROP:excep}. 
\end{proof} 

\ENDINPUT

%% file: Tables/Tables31.tex
\input hat
\NNULL{ }




\begin{table*}[htb]{}
{\scriptsize
$ rank(\Omega)= 2 $
$$
\def\arraystretch{1.4}
\begin{array}{|r|r|r|r|lr|lr|lr|lr|}
\hline
 d& c& v& w&E_6/T^2\bullet& f&E_7/T^2\bullet& f&E_8/T^2\bullet& f&
F_4/T^2\bullet& f
\\  \hline
14&144&4&36&{}&{}&{}&{}&{{A_1A_1A_2A_2}\,\,\,{\scriptstyle \left[ \sM 0&0\endsM\right]}}&2&{}&{}
\\  \hline
12&48&4&8&{}&{}&{}&{}&{{A_1A_2A_3}\,\,\,{\scriptstyle \left[ \sM 1&0\endsM\right]}}&4&{}&{}
\\  \hline
11&8&4&4&{}&{}&{}&{}&{{A_1A_1A_4}\,\,\,{\scriptstyle \left[ \sM 0&0\endsM\right]}}&4&{}&{}
\\  \hline
10&6&4&3&{}&{}&{}&{}&{{A_2A_4}\,\,\,{\scriptstyle \left[ \sM 1&1\endsM\right]}}&4&{}&{}
\\  \hline
10&4&4&1&{}&{}&{}&{}&{{A_3A_3}\,\,\,{\scriptstyle \left[ \sM 0&1\endsM\right]}}&2&{}&{}
\\  \hline
9&8&3&8&{}&{}&{{A_1A_1A_1A_2}\,\,\,{\scriptstyle \left[ \sM 0&0\endsM\right]}}&1&{{A_2D_4}{
 \left[ \sM 0&0\endsM\right]}}&1&{A_2}&1
\\  \hline
9&12&4&6&{}&{}&{}&{}&{{A_1A_5}\,\,\,{\scriptstyle \left[ \sM 1&1\endsM\right]}}&3&{}&{}
\\  \hline
8&4&3&2&{}&{}&{{A_1A_1A_3}\,\,\,{\scriptstyle \left[ \sM 0&0\endsM\right]}}&3&{{A_1D_5}\,\,\,{\scriptstyle \left[
 \sM 1&0\endsM\right]}}&3&{A_1\widetilde A_1}&3
\\  \hline
8&6&3&3&{}&{}&{{A_1A_2A_2}\,\,\,{\scriptstyle \left[ \sM 1&0\endsM\right]}}&3&{}&{}&{}&{}
\\  \hline
7&2&4&2&{}&{}&{{A_2A_3}\,\,\,{\scriptstyle \left[ \sM 0&1\endsM\right]}}&3&{{A_6}\,\,\,{\scriptstyle \left[ \sM 1&
1\endsM\right]}}&3&{}&{}
\\  \hline
6&1&3&1&{A_2A_2}&1&{{A_5''}\,\,\,{\scriptstyle \left[ \sM 0&0\endsM\right]}}&1&{{E_6}\,\,\,{\scriptstyle \left[
 \sM 1&0\endsM\right]}}&1&{\widetilde A_2}&1
\\  \hline
6&2&3&2&{A_1A_1A_2}&5&{{A_1A_4}\,\,\,{\scriptstyle \left[ \sM 1&1\endsM\right]}}&5&{}&{}&{}&{}
\\  \hline
6&4&2&4&{}&{}&{{A_1D_4}\,\,\,{\scriptstyle \left[ \sM 0&0\endsM\right]}}&1&{{D_6}\,\,\,{\scriptstyle \left[ \sM 1&
0\endsM\right]}}&1&{B_2}&1
\\  \hline
5&2&3&2&{A_1A_3}&4&{{A_5'}\,\,\,{\scriptstyle \left[ \sM 1&1\endsM\right]}}&2&{}&{}&{}&{}
\\  \hline
4&1&2&1&{A_4}&4&{{D_5}\,\,\,{\scriptstyle \left[ \sM 1&0\endsM\right]}}&2&{}&{}&{}&{}
\\  \hline
3&1&3&1&{D_4}&1&{}&{}&{}&{}&{}&{}
\\  \hline
\end{array}
$$
Note. We have $g_{\omega } \in \{1,2,3 \}$ for all $\omega\in \Omega$, except for two $T$-roots $+\omega, -\omega$ with $g_{\pm\omega }=4$ in the case $d=12$, $c=48=2^2\cdot 3 \cdot 4$.
}
\end{table*}
\begin{table*}[htb]{}
{\scriptsize
$ rank(\Omega)= 3 $
$$
\def\arraystretch{1.4}
\begin{array}{|r|r|r|r|lr|lr|lr|lr|}
\hline
 d& c& v& w&E_6/T^3\bullet& f&E_7/T^3\bullet& f&E_8/T^3\bullet& f&
F_4/T^3\bullet& f
\\  \hline
23&16&10&16&{}&{}&{}&{}&{A_1A_1A_1A_2\,\,\,{\scriptstyle \left[ \sM 1&0\endsM\right]}}&8&{}&{}
\\  \hline
21&6&10&3&{}&{}&{}&{}&{A_1A_2A_2\,\,\,{\scriptstyle \left[ \sM 1&1\endsM\right]}}&8&{}&{}
\\  \hline
20&8&10&4&{}&{}&{}&{}&{A_1A_1A_3\,\,\,{\scriptstyle \left[ \sM 1&1\endsM\right]}}&10&{}&{}
\\  \hline
18&2&11&2&{}&{}&{}&{}&{A_2A_3\,\,\,{\scriptstyle \left[ \sM 1&1\endsM\right]}}&10&{}&{}
\\  \hline
17&2&10&2&{}&{}&{}&{}&{A_1A_4\,\,\,{\scriptstyle \left[ \sM 1&1\endsM\right]}}&12&{}&{}
\\  \hline
16&8&7&8&{}&{}&{4A_1\,\,\,{\scriptstyle \left[ \sM 0&0\endsM\right]}}&2&{A_1D_4\,\,\,{\scriptstyle \left[ \sM 1&0
\endsM\right]}}&2&A_1&2
\\  \hline
14&2&10&2&{}&{}&{}&{}&{A_5\vphantom{footnote}\,\,\,{\scriptstyle \left[ \sM 1&1\endsM\right]}}&4
&{}&{}
\\  \hline
14&2&8&2&{}&{}&{A_1A_1A_2\,\,\,{\scriptstyle \left[ \sM 1&1\endsM\right]}}&12&{}&{}&{}&{}
\\  \hline
13&1&9&1&{}&{}&{A_2A_2\,\,\,{\scriptstyle \left[ \sM 1&1\endsM\right]}}&4&{}&{}&{}&{}
\\  \hline
13&1&6&1&{}&{}&{[A_1A_3]''\,\,\,{\scriptstyle \left[ \sM 0&0\endsM\right]}}&2&{D_5\,\,\,{\scriptstyle \left[ \sM 1
&0\endsM\right]}}&2&\widetilde A_1&2
\\  \hline
12&2&8&2&{}&{}&{[A_1A_3]'\,\,\,{\scriptstyle \left[ \sM 1&1\endsM\right]}}&9&{}&{}&{}&{}
\\  \hline
11&2&7&2&A_1A_1A_1&5&{}&{}&{}&{}&{}&{}
\\  \hline
10&1&7&1&A_1A_2&10&{A_4\,\,\,{\scriptstyle \left[ \sM 1&1\endsM\right]}}&5&{}&{}&{}&{}
\\  \hline
9&1&3&1&{}&{}&{D_4\,\,\,{\scriptstyle \left[ \sM 1&0\endsM\right]}}&1&{}&{}&{}&{}
\\  \hline
8&1&6&1&A_3&5&{}&{}&{}&{}&{}&{}
\\  \hline
\end{array}
$$
Note. We have $g_{\omega } \in \{1,2\}$ for all $\omega\in \Omega$, except for two $T$-roots $+\omega, -\omega$ with $g_{\pm\omega }=3$ in the case $d=21$, $c=6$.
}
\end{table*}
\begin{table*}[htb]{}
{\scriptsize
$ rank(\Omega)= 4 $
$$
\def\arraystretch{1.4}
\begin{array}{|r|r|r|r|lr|lr|lr|lr|}
\hline
 d& c& v& w&E_6/T^4\bullet& f&E_7/T^4\bullet& f&E_8/T^4\bullet& f&
F_4/T^4\bullet& f
\\  \hline
36&16&20&16&{}&{}&{}&{}&{4A_1\,\,\,{\scriptstyle \left[ \sM 1&1\endsM\right]}}&7&{}&{}
\\  \hline
33&2&21&2&{}&{}&{}&{}&{A_1A_1A_2\,\,\,{\scriptstyle \left[ \sM 1&1\endsM\right]}}&28&{}&{}
\\  \hline
30&1&24&1&{}&{}&{}&{}&{A_2A_2\,\,\,{\scriptstyle \left[ \sM 1&1\endsM\right]}}&8&{}&{}
\\  \hline
29&2&21&2&{}&{}&{}&{}&{A_1A_3\,\,\,{\scriptstyle \left[ \sM 1&1\endsM\right]}}&20&{}&{}
\\  \hline
25&1&20&1&{}&{}&{}&{}&{A_4\,\,\,{\scriptstyle \left[ \sM 1&1\endsM\right]}}&6&{}&{}
\\  \hline
24&1&12&1&{}&{}&{A_1A_1A_1\,\,\,{\scriptstyle \left[ \sM 0&0\endsM\right]}}&1&{D_4\,\,\,{\scriptstyle \left[ \sM 1
&0\endsM\right]}}&1&\{e\}&1
\\  \hline
23&2&16&2&{}&{}&{A_1A_1A_1\,\,\,{\scriptstyle \left[ \sM 1&1\endsM\right]}}&10&{}&{}&{}&{}
\\  \hline
21&1&17&1&{}&{}&{A_1A_2\,\,\,{\scriptstyle \left[ \sM 1&1\endsM\right]}}&18&{}&{}&{}&{}
\\  \hline
18&1&15&1&{}&{}&{A_3\,\,\,{\scriptstyle \left[ \sM 1&1\endsM\right]}}&6&{}&{}&{}&{}
\\  \hline
17&1&14&1&A_1A_1&10&{}&{}&{}&{}&{}&{}
\\  \hline
15&1&15&1&A_2&5&{}&{}&{}&{}&{}&{}
\\  \hline
\end{array}
$$
}
\end{table*}
\begin{table*}[htb]{}
{\scriptsize
$ rank(\Omega)= 5 $
$$
\def\arraystretch{1.4}
\begin{array}{|r|r|r|r|lr|lr|lr|}
\hline
 d& c& v& w&E_6/T^5\bullet& f&E_7/T^5\bullet& f&E_8/T^5\bullet& f
\\  \hline
50&2&37&2&{}&{}&{}&{}&{A_1A_1A_1\,\,\,{\scriptstyle \left[ \sM 1&1\endsM\right]}}&21
\\  \hline
46&1&40&1&{}&{}&{}&{}&{A_1A_2\,\,\,{\scriptstyle \left[ \sM 1&1\endsM\right]}}&28
\\  \hline
41&1&36&1&{}&{}&{}&{}&{A_3\,\,\,{\scriptstyle \left[ \sM 1&1\endsM\right]}}&7
\\  \hline
33&1&29&1&{}&{}&{A_1A_1\,\,\,{\scriptstyle \left[ \sM 1&1\endsM\right]}}&15&{}&{}
\\  \hline
30&1&30&1&{}&{}&{A_2\,\,\,{\scriptstyle \left[ \sM 1&1\endsM\right]}}&6&{}&{}
\\  \hline
25&1&25&1&A_1&6&{}&{}&{}&{}
\\  \hline
\end{array}
$$
}
\end{table*}
\begin{table*}[htb]{}
{\scriptsize
$ rank(\Omega)= 6 $
$$
\def\arraystretch{1.4}
\begin{array}{|r|r|r|r|lr|lr|lr|}
\hline
 d& c& v& w&E_6/T^6\bullet& f&E_7/T^6\bullet& f&E_8/T^6\bullet& f
\\  \hline
68&1&62&1&{}&{}&{}&{}&{A_1A_1,\,\,{\scriptstyle \left[ \sM 1&1\endsM\right]}}&21
\\  \hline
63&1&63&1&{}&{}&{}&{}&{A_2\,\,\,{\scriptstyle \left[ \sM 1&1\endsM\right]}}&7
\\  \hline
46&1&46&1&{}&{}&{A_1\,\,\,{\scriptstyle \left[ \sM 1&1\endsM\right]}}&7&{}&{}
\\  \hline
36&1&36&1&\{e\}&1&{}&{}&{}&{}
\\  \hline
\end{array}
$$
}
{\scriptsize
$ rank(\Omega)= 7 $
$$
\def\arraystretch{1.4}
\begin{array}{|r|r|r|r|lr|lr|lr|}
\hline
 d& c& v& w&E_6/T^7\bullet& f&E_7/T^7\bullet& f&E_8/T^7\bullet& f
\\  \hline
91&1&91&1&{}&{}&{}&{}&{A_1\,\,\,{\scriptstyle \left[ \sM 1&1\endsM\right]}}&8
\\  \hline
63&1&63&1&{}&{}&{\{e\}\,\,\,{\scriptstyle \left[ \sM 1&1\endsM\right]}}&1&{}&{}
\\  \hline
\end{array}
$$
}
{\scriptsize
$ rank(\Omega)= 8 $
$$
\def\arraystretch{1.4}
\begin{array}{|r|r|r|r|lr|lr|lr|}
\hline
 d& c& v& w&E_6/T^8\bullet& f&E_7/T^8\bullet& f&E_8/T^8\bullet& f
\\  \hline
120&1&120&1&{}&{}&{}&{}&{\{e\}\,\,\,{\scriptstyle \left[ \sM 1&1\endsM\right]}}&1
\\  \hline
\end{array}
$$
}
\end{table*}


%% file: Tables/Nmbering.tex

\begin{table*}[htb]{}
\caption{Consistent numberings of simple T-roots}
{\scriptsize
\begin{quote}
Consistent numberings of simple T-roots (in each line).
We skip $5$ of $10$ diagrams of the type $E_6$. These $5$ diagrams may be obtained by symmetry of the Dynkin graph.
\end{quote}
$$
\label{table:Numberings}
\def\arraystretch{1.4}
\begin{array}{|r|r|r|r|r|lr|lr|lr|lr|}
\hline
k& d& c& v& w&E_6/T^kH& f&E_7/T^kH& f&E_8/T^kH& f&F_4/T^kH& f
\\  \hline
4&24&1&12&1&{}&{}&{ 
\begin{array}{l}E_{7}/T^4[A_1A_1A_1]''\\
{\sM 0&1&0&2&3&4\\ &&&0\endsM\ig}
\end{array}}&1&{
\begin{array}{l}E_{8}/T^4D_4\\
{\sM 4&3&2&0&0&0&1\\ &&&&0\endsM\ig}
\end{array}}&1&{
\begin{array}{l}F_{4}/T^4\\
{\sM 1&2\Leftarrow3&4\endsM\ig}
\end{array}}&1
\\  \hline
3&16&8&7&8&{}&{}&{ \begin{array}{l}E_{7}/T^3A_1A_1A_1A_1\\{\sM 0&1&0&2&3&
0\\ &&&0\endsM\ig}\\{\sM 0&1&0&2&0&3\\ &&&0\endsM\ig}\end{array}}&2&{
 \begin{array}{l}E_{8}/T^3A_1D_4\\{\sM 0&3&2&0&0&0&1\\ &&&&0\endsM\ig}\\{\sM
3&0&2&0&0&0&1\\ &&&&0\endsM\ig}\end{array}}&2&{ \begin{array}{l}F_4/T^3A_1\\{
\sM 1&2\Leftarrow3&0\endsM\ig}\\{\sM 1&2\Leftarrow0&3\endsM\ig}\end{array}}&2
\\  \hline
3&13&1&6&1&{}&{}&{ \begin{array}{l}E_{7}/T^3[A_1A_3]''\\{\sM 0&1&0&0&2&3\\ &&&
0\endsM\ig}\\{\sM 0&0&0&1&2&3\\ &&&0\endsM\ig}\end{array}}&2&{
 \begin{array}{l}E_{8}/T^3D_5\\{\sM 3&2&0&0&0&0&1\\ &&&&0\endsM\ig}\\{\sM 3&2&
1&0&0&0&0\\ &&&&0\endsM\ig}\end{array}}&2&{ \begin{array}{l}
F_4/T^3\widetilde A_1\\{\sM 1&0\Leftarrow2&3\endsM\ig}\\{\sM 0&1\Leftarrow2&3
\endsM\ig}\end{array}}&2
\\  \hline
3&10&1&7&1&{ \begin{array}{l}E_{6}/T^3A_1A_2\!\!\!\!\!\!\!\!\! \\
{\sM 1&2&3&0&0\\ &&0\endsM\ig}\\ 
{\sM 1&2&0&3&0\\ &&0\endsM\ig}\\
{\sM 1&0&2&0&0\\ &&3\endsM\ig}\\
{\sM 1&0&0&2&0\\ &&3\endsM\ig}\\
   {\sM 0&0&3&2&0\\ &&1\endsM\ig}
\end{array}}&10&{
 \begin{array}{l}E_{7}/T^3A_4\\
 {\sM 1&2&3&0&0&0\\ &&&0\endsM\ig}\\
 {\sM 1&2&0&0&0&0\\ &&&3\endsM\ig}\\ 
 {\sM 1&0&0&0&2&3\\ &&&0\endsM\ig}\\
 {\sM 1&0&0&0&0&3\\ &&&2\endsM\ig}\\
 {\sM 0&0&0&0&3&1\\ &&&2\endsM\ig}
 \end{array}}&5&{}&{}&{}&{}
\\  \hline
2&8&4&3&2&{}&{}&{ \begin{array}{l}E_{7}/T^2{A_1A_1A_3}\\{\sM 0&1&0&0&2&0\\ &&&
0\endsM\ig}\\{\sM 0&0&0&1&2&0\\ &&&0\endsM\ig}\\{\sM 0&0&0&1&0&2\\ &&&0
\endsM\ig}\end{array}}&3&{ \begin{array}{l}E_{8}/T^2{A_1D_5}\\{\sM 0&2&0&0&0&
0&1\\ &&&&0\endsM\ig}\\{\sM 0&2&1&0&0&0&0\\ &&&&0\endsM\ig}\\{\sM 2&0&1&0&0&0&
0\\ &&&&0\endsM\ig}\end{array}}&3&{ \begin{array}{l}F_4/T^2{A_1\widetilde A_1}\!\!\!\!\!\!\!
\\{\sM 1&0\Leftarrow2&0\endsM\ig}\\{\sM 0&1\Leftarrow2&0\endsM\ig}\\{\sM 0&1
\Leftarrow0&2\endsM\ig}\end{array}}&3
\\  \hline
2&7&2&4&2&{}&{}&{ \begin{array}{l}E_{7}/T^2{A_2A_3}\\{\sM 0&0&1&0&0&2\\ &&&0
\endsM\ig}\\{\sM 0&0&1&0&0&0\\ &&&2\endsM\ig}\\{\sM 0&0&0&1&0&0\\ &&&2
\endsM\ig}\end{array}}&3&{ \begin{array}{l}E_{8}/T^2{A_6}\\{\sM 2&0&0&0&0&0&
0\\ &&&&1\endsM\ig}\\{\sM 0&0&0&0&0&0&2\\ &&&&1\endsM\ig}\\{\sM 0&0&0&0&0&1&
2\\ &&&&0\endsM\ig}\end{array}}&3&{}&{}
\\  \hline
\end{array}
$$
}
\end{table*}

%% file: Tables/Tables32.tex
\input hat
\NNULL{ }






\begin{table*}[htb]{}
\caption{The same table as 4 with diagrams of some extremal complex forms} 

\begin{quote}
{\scriptsize The same, as Table~\ref{table:31}, with diagrams of some standard complex forms of flag manifolds. The positive labels (on each diagram) are coefficients of the major $T$-root. We skip only the diagrams with coprime labels. For $rank(\Omega )>5$ all diagrams are skipped.}%
\end{quote}

\label{table:32}

\vskip2ex

{\scriptsize
$ rank(\Omega)= 2 $
$$
\def\arraystretch{1.4}
\begin{array}{|r|r|r|r|lr|lr|lr|lr|}
\hline
 d& c& v& w&E_6/T^2\bullet& f&E_7/T^2\bullet& f&E_8/T^2\bullet& f&
F_4/T^2\bullet& f
\\  \hline
14&144&4&36&{}&{}&{}&{}&{ \begin{array}{l}{A_1A_1A_2A_2}\\{\sM 0&3&0&0&6&0&
0\\ &&&&0\endsM\ig}\\{\sM 0&0&4&0&6&0&0\\ &&&&0\endsM\ig}\end{array}}&2&{}&{}
\\  \hline
12&48&4&8&{}&{}&{}&{}&{ \begin{array}{l}{A_1A_2A_3}\\{\sM 2&0&0&0&6&0&0\\ &&&&
0\endsM\ig}\\{\sM 0&0&4&0&0&4&0\\ &&&&0\endsM\ig}\end{array}}&4&{}&{}
\\  \hline
11&8&4&4&{}&{}&{}&{}&{ \begin{array}{l}{A_1A_1A_4}\\{\sM 0&0&0&0&6&4&0\\ &&&&0
\endsM\ig}\\{\sM 0&0&0&0&6&0&2\\ &&&&0\endsM\ig}\end{array}}&4&{}&{}
\\  \hline
10&6&4&3&{}&{}&{}&{}&{ \begin{array}{l}{A_2A_4}\\{\sM 0&0&0&0&6&0&0\\ &&&&3
\endsM\ig}\end{array}}&4&{}&{}
\\  \hline
10&4&4&1&{}&{}&{}&{}&{ \begin{array}{l}{A_3A_3}\end{array}}&2&{}&{}
\\  \hline
9&8&3&8&{}&{}&{ \begin{array}{l}{A_1A_1A_1A_2}\\{\sM 0&2&0&4&0&0\\ &&&0
\endsM\ig}\end{array}}&1&{ \begin{array}{l}{A_2D_4}\\{\sM 0&0&4&0&0&0&2\\ &&&&
0\endsM\ig}\end{array}}&1&{ \begin{array}{l}{A_2}\\{\sM 2&4\Leftarrow0&0
\endsM\ig}\end{array}}&1
\\  \hline
9&12&4&6&{}&{}&{}&{}&{ \begin{array}{l}{A_1A_5}\\{\sM 2&0&0&0&0&4&0\\ &&&&0
\endsM\ig}\\{\sM 0&3&0&0&0&0&0\\ &&&&3\endsM\ig}\end{array}}&3&{}&{}
\\  \hline
8&4&3&2&{}&{}&{ \begin{array}{l}{A_1A_1A_3}\\{\sM 0&0&0&4&0&2\\ &&&0\endsM\ig}
\end{array}}&3&{ \begin{array}{l}{A_1D_5}\\{\sM 2&0&4&0&0&0&0\\ &&&&0\endsM\ig
}\end{array}}&3&{ \begin{array}{l}{A_1\widetilde A_1}\\{\sM 0&4\Leftarrow0&2
\endsM\ig}\end{array}}&3
\\  \hline
8&6&3&3&{}&{}&{ \begin{array}{l}{A_1A_2A_2}\\{\sM 0&0&3&0&3&0\\ &&&0\endsM\ig}
\end{array}}&3&{}&{}&{}&{}
\\  \hline
7&2&4&2&{}&{}&{ \begin{array}{l}{A_2A_3}\\{\sM 0&0&0&4&0&0\\ &&&2\endsM\ig}
\end{array}}&3&{ \begin{array}{l}{A_6}\\{\sM 0&0&0&0&0&4&2\\ &&&&0\endsM\ig}
\end{array}}&3&{}&{}
\\  \hline
6&2&3&2&{ \begin{array}{l}{A_1A_1A_2}\\{\sM 0&2&0&2&0\\ &&0\endsM\ig}
\end{array}}&5&{ \begin{array}{l}{A_1A_4}\\{\sM 0&2&0&0&0&0\\ &&&2\endsM\ig}
\end{array}}&5&{}&{}&{}&{}
\\  \hline
6&4&2&4&{}&{}&{ \begin{array}{l}{A_1D_4}\\{\sM 0&2&0&0&0&2\\ &&&0\endsM\ig}
\end{array}}&1&{ \begin{array}{l}{D_6}\\{\sM 2&0&0&0&0&0&2\\ &&&&0\endsM\ig}
\end{array}}&1&{ \begin{array}{l}{B_2}\\{\sM 2&0\Leftarrow0&2\endsM\ig}
\end{array}}&1
\\  \hline
6&1&3&1&{ \begin{array}{l}{A_2A_2}\end{array}}&1&{ \begin{array}{l}{A_5''}
\end{array}}&1&{ \begin{array}{l}{E_6}\end{array}}&1&{ \begin{array}{l}
{\widetilde A_2}\end{array}}&1
\\  \hline
5&2&3&2&{ \begin{array}{l}{A_1A_3}\\{\sM 0&2&0&0&0\\ &&2\endsM\ig}\\{\sM 0&0&
0&2&0\\ &&2\endsM\ig}\end{array}}&4&{ \begin{array}{l}{A_5'}\\{\sM 0&0&0&0&0&
2\\ &&&2\endsM\ig}\end{array}}&2&{}&{}&{}&{}
\\  \hline
4&1&2&1&{ \begin{array}{l}{A_4}\end{array}}&4&{ \begin{array}{l}{D_5}
\end{array}}&2&{}&{}&{}&{}
\\  \hline
3&1&3&1&{ \begin{array}{l}{D_4}\end{array}}&1&{}&{}&{}&{}&{}&{}
\\  \hline
\end{array}
$$
}
\end{table*}
\begin{table*}[htb]{}
{\scriptsize
$ rank(\Omega)= 3 $
$$
\def\arraystretch{1.4}
\begin{array}{|r|r|r|r|lr|lr|lr|lr|}
\hline
 d& c& v& w&E_6/T^3\bullet& f&E_7/T^3\bullet& f&E_8/T^3\bullet& f&
F_4/T^3\bullet& f
\\  \hline
23&16&10&16&{}&{}&{}&{}&{ \begin{array}{l}A_1A_1A_1A_2\\{\sM 2&0&4&0&6&0&
0\\ &&&&0\endsM\ig}\\{\sM 0&0&4&0&6&4&0\\ &&&&0\endsM\ig}\\{\sM 0&0&4&0&6&0&
2\\ &&&&0\endsM\ig}\end{array}}&8&{}&{}
\\  \hline
21&6&10&3&{}&{}&{}&{}&{ \begin{array}{l}A_1A_2A_2\\{\sM 0&3&0&0&6&0&0\\ &&&&3
\endsM\ig}\end{array}}&8&{}&{}
\\  \hline
20&8&10&4&{}&{}&{}&{}&{ \begin{array}{l}A_1A_1A_3\\{\sM 2&0&4&0&0&4&0\\ &&&&0
\endsM\ig}\\{\sM 2&0&0&0&6&4&0\\ &&&&0\endsM\ig}\\{\sM 2&0&0&0&6&0&2\\ &&&&0
\endsM\ig}\end{array}}&10&{}&{}
\\  \hline
18&2&11&2&{}&{}&{}&{}&{ \begin{array}{l}A_2A_3\\{\sM 0&0&4&0&0&4&2\\ &&&&0
\endsM\ig}\end{array}}&10&{}&{}
\\  \hline
17&2&10&2&{}&{}&{}&{}&{ \begin{array}{l}A_1A_4\\{\sM 0&0&0&0&6&4&2\\ &&&&0
\endsM\ig}\end{array}}&12&{}&{}
\\  \hline
16&8&7&8&{}&{}&{ \begin{array}{l}A_1A_1A_1A_1\\{\sM 0&2&0&4&0&2\\ &&&0
\endsM\ig}\end{array}}&2&{ \begin{array}{l}A_1D_4\\{\sM 2&0&4&0&0&0&2\\ &&&&0
\endsM\ig}\end{array}}&2&{ \begin{array}{l}A_1\\{\sM 2&4\Leftarrow0&2\endsM\ig
}\end{array}}&2
\\  \hline
14&2&8&2&{}&{}&{ \begin{array}{l}A_1A_1A_2\\{\sM 0&2&0&4&0&0\\ &&&2\endsM\ig}
\end{array}}&12&{}&{}&{}&{}
\\  \hline
14&2&10&2&{}&{}&{}&{}&{ \begin{array}{l}A_5\vphantom{footnote}\\{\sM 2&0&0&0&
0&4&2\\ &&&&0\endsM\ig}\end{array}}&4&{}&{}
\\  \hline
13&1&9&1&{}&{}&{ \begin{array}{l}A_2A_2\end{array}}&4&{}&{}&{}&{}
\\  \hline
13&1&6&1&{}&{}&{ \begin{array}{l}[A_1A_3]''\end{array}}&2&{ \begin{array}{l}
D_5\end{array}}&2&{ \begin{array}{l}\widetilde A_1\end{array}}&2
\\  \hline
12&2&8&2&{}&{}&{ \begin{array}{l}[A_1A_3]'\\{\sM 0&2&0&0&0&2\\ &&&2\endsM\ig}
\\{\sM 0&0&0&4&0&2\\ &&&2\endsM\ig}\end{array}}&9&{}&{}&{}&{}
\\  \hline
11&2&7&2&{ \begin{array}{l}A_1A_1A_1\\{\sM 0&2&0&2&0\\ &&2\endsM\ig}
\end{array}}&5&{}&{}&{}&{}&{}&{}
\\  \hline
10&1&7&1&{ \begin{array}{l}A_1A_2\end{array}}&10&{ \begin{array}{l}A_4
\end{array}}&5&{}&{}&{}&{}
\\  \hline
9&1&3&1&{}&{}&{ \begin{array}{l}D_4\end{array}}&1&{}&{}&{}&{}
\\  \hline
8&1&6&1&{ \begin{array}{l}A_3\end{array}}&5&{}&{}&{}&{}&{}&{}
\\  \hline
\end{array}
$$
}
\end{table*}
\begin{table*}[htb]{}
{\scriptsize
$ rank(\Omega)= 4 $
$$
\def\arraystretch{1.4}
\begin{array}{|r|r|r|r|lr|lr|lr|lr|}
\hline
 d& c& v& w&E_6/T^4\bullet& f&E_7/T^4\bullet& f&E_8/T^4\bullet& f&
F_4/T^4\bullet& f
\\  \hline
36&16&20&16&{}&{}&{}&{}&{ \begin{array}{l}A_1A_1A_1A_1\\{\sM 2&0&4&0&6&4&
0\\ &&&&0\endsM\ig}\\{\sM 2&0&4&0&6&0&2\\ &&&&0\endsM\ig}\end{array}}&7&{}&{}
\\  \hline
33&2&21&2&{}&{}&{}&{}&{ \begin{array}{l}A_1A_1A_2\\{\sM 0&0&4&0&6&4&2\\ &&&&0
\endsM\ig}\end{array}}&28&{}&{}
\\  \hline
30&1&24&1&{}&{}&{}&{}&{ \begin{array}{l}A_2A_2\end{array}}&8&{}&{}
\\  \hline
29&2&21&2&{}&{}&{}&{}&{ \begin{array}{l}A_1A_3\\{\sM 2&0&4&0&0&4&2\\ &&&&0
\endsM\ig}\\{\sM 2&0&0&0&6&4&2\\ &&&&0\endsM\ig}\end{array}}&20&{}&{}
\\  \hline
25&1&20&1&{}&{}&{}&{}&{ \begin{array}{l}A_4\end{array}}&6&{}&{}
\\  \hline
24&1&12&1&{}&{}&{ \begin{array}{l}[A_1A_1A_1]''\end{array}}&1&{ \begin{array}{l}
D_4\end{array}}&1&{ \begin{array}{l}\{e\}\end{array}}&1
\\  \hline
23&2&16&2&{}&{}&{ \begin{array}{l}[A_1A_1A_1]'\\{\sM 0&2&0&4&0&2\\ &&&2\endsM\ig}
\end{array}}&10&{}&{}&{}&{}
\\  \hline
21&1&17&1&{}&{}&{ \begin{array}{l}A_1A_2\end{array}}&18&{}&{}&{}&{}
\\  \hline
18&1&15&1&{}&{}&{ \begin{array}{l}A_3\end{array}}&6&{}&{}&{}&{}
\\  \hline
17&1&14&1&{ \begin{array}{l}A_1A_1\end{array}}&10&{}&{}&{}&{}&{}&{}
\\  \hline
15&1&15&1&{ \begin{array}{l}A_2\end{array}}&5&{}&{}&{}&{}&{}&{}
\\  \hline
\end{array}
$$
}
\end{table*}
\begin{table*}[htb]{}
{\scriptsize
$ rank(\Omega)= 5 $
$$
\def\arraystretch{1.4}
\begin{array}{|r|r|r|r|lr|lr|lr|}
\hline
 d& c& v& w&E_6/T^5\bullet& f&E_7/T^5\bullet& f&E_8/T^5\bullet& f
\\  \hline
50&2&37&2&{}&{}&{}&{}&{ \begin{array}{l}A_1A_1A_1\\{\sM 2&0&4&0&6&4&2\\ &&&&0
\endsM\ig}\end{array}}&21
\\  \hline
46&1&40&1&{}&{}&{}&{}&{ \begin{array}{l}A_1A_2\end{array}}&28
\\  \hline
41&1&36&1&{}&{}&{}&{}&{ \begin{array}{l}A_3\end{array}}&7
\\  \hline
33&1&29&1&{}&{}&{ \begin{array}{l}A_1A_1\end{array}}&15&{}&{}
\\  \hline
30&1&30&1&{}&{}&{ \begin{array}{l}A_2\end{array}}&6&{}&{}
\\  \hline
25&1&25&1&{ \begin{array}{l}A_1\end{array}}&6&{}&{}&{}&{}
\\  \hline
\end{array}
$$
}
\end{table*}
\begin{table*}[htb]{}
{\scriptsize
$ rank(\Omega)= 6 $
$$
\def\arraystretch{1.4}
\begin{array}{|r|r|r|r|lr|lr|lr|}
\hline
 d& c& v& w&E_6/T^6\bullet& f&E_7/T^6\bullet& f&E_8/T^6\bullet& f
\\  \hline
68&1&62&1&{}&{}&{}&{}&{ \begin{array}{l}A_1A_1\end{array}}&21
\\  \hline
63&1&63&1&{}&{}&{}&{}&{ \begin{array}{l}A_2\end{array}}&7
\\  \hline
46&1&46&1&{}&{}&{ \begin{array}{l}A_1\end{array}}&7&{}&{}
\\  \hline
36&1&36&1&{ \begin{array}{l}\{e\}\end{array}}&1&{}&{}&{}&{}
\\  \hline
\end{array}
$$
}
{\scriptsize
$ rank(\Omega)= 7 $
$$
\def\arraystretch{1.4}
\begin{array}{|r|r|r|r|lr|lr|lr|}
\hline
 d& c& v& w&E_6/T^7\bullet& f&E_7/T^7\bullet& f&E_8/T^7\bullet& f
\\  \hline
91&1&91&1&{}&{}&{}&{}&{ \begin{array}{l}A_1\end{array}}&8
\\  \hline
63&1&63&1&{}&{}&{ \begin{array}{l}\{e\}\end{array}}&1&{}&{}
\\  \hline
\end{array}
$$
}
{\scriptsize
$ rank(\Omega)= 8 $
$$
\def\arraystretch{1.4}
\begin{array}{|r|r|r|r|lr|lr|lr|}
\hline
 d& c& v& w&E_6/T^8\bullet& f&E_7/T^8\bullet& f&E_8/T^8\bullet& f
\\  \hline
120&1&120&1&{}&{}&{}&{}&{ \begin{array}{l}\{e\}\end{array}}&1
\\  \hline
\end{array}
$$
}
\end{table*}

\ENDINPUT